\documentclass[11pt,letter,final]{amsart}
\newcommand{\withPgf}{no}
\newcommand{\atype}{preprint}

\usepackage{xr}
\usepackage{hyperref}
\usepackage{booktabs,colortbl}
\usepackage{framed, color}
\usepackage{amsmath}
\usepackage{bm}
\usepackage{amsfonts}
\usepackage{color}
\usepackage{algorithm}
\usepackage{algpseudocode}

\ifthenelse{\equal{\atype}{siam}}{
\bibliographystyle{siam}
}{%
\usepackage{amsthm}
\bibliographystyle{siam}

}

\newtheorem{remark}{Remark}[section]

\ifthenelse{\equal{\atype}{siam}}
{
\newcommand{\itemname}{romannum}
}
{

\newcommand{\itemname}{itemize}
}

\definecolor{Blue}{rgb}{0,0,1}
\definecolor{Red}{rgb}{1,0,0}
\definecolor{Orange}{rgb}{1,0,1}

\newcommand{\cb}[1]{\left\{#1\right\}}
\newcommand{\rb}[1]{\left(#1\right)}
\newcommand{\eb}[1]{\left[#1\right]}
\newcommand{\vertiii}[1]{{\left\vert\kern-0.25ex\left\vert\kern-0.25ex\left\vert #1 
    \right\vert\kern-0.25ex\right\vert\kern-0.25ex\right\vert}}
\newcommand{\detI}{\ensuremath{\text{h}}}
\newcommand{\redI}{\ensuremath{\text{red}}}
\newcommand{\norm}[1]{\left\|#1\right\|}
\newcommand{\enorm}[1]{\vertiii{#1}}
\newcommand{\xnorm}[1]{\norm{#1}_{X_{\detI}}}
\newcommand{\K}{\mathbf{K}}

\newcommand{\by}{\mathbf{y}}
\newcommand{\bI}{\mathbf{I}}
\newcommand{\bx}{\mathbf{x}}
\newcommand{\bw}{\mathbf{w}}
\newcommand{\abs}[1]{\left|#1\right|}

\def \cal{\mathcal}

\newcommand{\bmu}{{\bm \mu}}
\newcommand{\inputs}{{\bmu}}
\newcommand{\params}{{\inputs}}
\newcommand{\paramTrain}[1]{{\inputs}_{#1}}

\newcommand{\nTrain}{N}
\newcommand{\bmutrain}{\tilde{\bm \mu}}

\newcommand{\field}[1]{ \mathbb{#1} }
\usepackage{mathtools}
\newcommand{\defeq}{\vcentcolon=}
\newcommand{\eqdef}{=\vcentcolon}
\newcommand{\leftbasis}{\ensuremath{\bm{W}}}
\newcommand{\rightbasis}{\ensuremath{\bm{V}}} 
\newcommand{\rightbasisfine}{} 
\newcommand{\redstatefine}{\state} 
\newcommand{\rightbasiscoarse}{\rightbasis} 
\newcommand{\rightbasisdualcoarse}{\bm{Y}} 
 
\newcommand{\rightbasisdualcoarseT}{\rightbasisdualcoarse^T} 
\newcommand{\redstatecoarse}{\redstatecoords}

\newcommand{\rightbasisfineT}{} 
\newcommand{\range}[1]{\ensuremath{\mathrm{range}\left(#1\right)}}
\newcommand{\RR}{ \field{R} }

\newcommand{\cP}{ {\cal P} }
\newcommand{\cM}{ {\cal M} }

\newcommand{\cN}{ {\cal N} }

\newcommand{\cL}{ {\cal L} }

\newcommand{\myspan}{ \mathrm{span} }

\newcommand{\card}[1]{\mathrm{card}\left(#1\right)}

\newcommand{\erf}{\mathrm{erf}}
\newcommand{\rigorbound}{\quantity_\mathrm{LB}}
\newcommand{\crigor}{c}
\newcommand{\cvalidation}{c_\text{validation}}
\newcommand{\mode}[1]{\mathrm{mode}\rb{#1}}
\newcommand{\median}[1]{\mathrm{median}\rb{#1}}
\newcommand{\contstate}{u}
\newcommand{\contstatedummy}{w}
\newcommand{\contfestate}{\contstate_{\detI}}
\newcommand{\contfedual}{y_{\detI}}
\newcommand{\contfetest}{v_{\detI}}
\newcommand{\contredstate}{\contstate_{\redI}}
\newcommand{\state}{{\bm{u}}}

\newcommand{\redstate}{\state_{\redI}}
\newcommand{\dualstate}{{\bm{y}}}

\newcommand{\redstatecoords}{\hat \state}

\newcommand{\out}{s}
\newcommand{\outputs}{\out}
\newcommand{\outputfun}{g}
\newcommand{\outputfunCont}{\bar g}

\newcommand{\residual}{\ensuremath{\bm{r}}}
\newcommand{\modresidual}{\residual}
\newcommand{\coarsetofine}{\rightbasis}
\newcommand{\adjoint}{\bm{y}}
\newcommand{\redadjointfine}{\adjoint}
\newcommand{\redadjointcoarsetofine}{\adjoint_{\redI}}
\newcommand{\redadjointcoarsetofineII}{\adjoint_{\redI, i}}
\newcommand{\redadjointcoarsetofinearg}[1]{\adjoint_{\redI, #1}}

\newcommand{\redadjointcoarse}{\hat\adjoint}
\newcommand{\redadjointfineT}{\adjoint^T}
\newcommand{\nrbfine}{\ndof}

\newcommand{\residualnorm}{\ensuremath{{r}}}
\newcommand{\residualhyper}{\ensuremath{\tilde\residual}}
\newcommand{\outmeas}{\bar s}
\newcommand{\redout}{\out_{\redI}}

\newcommand{\outsurr}{\out_{\mathrm{surr}}}
\newcommand{\error}{\delta}
\newcommand{\quantity}{\ensuremath{m}}
\newcommand{\quantitysurr}{\ensuremath{\tilde{\quantity}}}

\newcommand{\stateerrorvecNo}{{\bm{\delta}_{\state}}}
\newcommand{\stateerrorvec}[1]{\stateerrorvecNo#1}
\newcommand{\stateerror}[1]{\norm{\stateerrorvec{#1}}}
\newcommand{\xstateerror}[1]{\xnorm{\stateerrorvec{#1}}}
\newcommand{\enstateerror}[1]{\enorm{\stateerrorvec{#1}}}
\newcommand{\outbias}{\delta_{\out}}
\newcommand{\outbiasone}{\delta_{\out_1}}
\newcommand{\outbiastwo}{\delta_{\out_2}}

\newcommand{\outbiasMF}{\delta_{\out, \text{MF}}}
\newcommand{\outerror}[1]{\abs{\delta_{\out}{#1}}}
\newcommand{\bound}{\Delta}
\newcommand{\lowerbound}{\Delta^\text{LB}}
\newcommand{\statebound}{\bound_{\state}}
\newcommand{\statelowerbound}{\lowerbound_{\state}}
\usepackage{bbm}
\newcommand{\indicator}[2]{{\mathbbm{1}}_{#2}\left(#1\right)}
\newcommand{\cint}{C}
\newcommand{\confidenceInterval}[2]{\cint_{#1}\left(#2\right)}
\newcommand{\confidenceIntervalOne}[1]{\cint_{#1}}
\newcommand{\indicators}{\bm{\rho}}
\newcommand{\stateindicator}{\bm{\rho}}
\newcommand{\outbound}{\bound_{\out}}
\newcommand{\outlowerbound}{\lowerbound_{\out}}
\newcommand{\outcorrected}{\out_{\text{corr}}}

\newcommand{\nindicator}{q}
\newcommand{\Para}{\cP}
\newcommand{\paramDomain}{\Para}
\newcommand{\Paratrain}{\ensuremath{\Para_{\mathrm{train}}}}
\newcommand{\Paragreedy}{\ensuremath{\Para_{\mathrm{greedy}}}}
\newcommand{\romesTraining}{\ensuremath{N}}
\newcommand{\Pararomes}{\ensuremath{\bar\Para}}
\newcommand{\Paratrainromes}{\ensuremath{\Para_{\mathrm{learn}}}}
\newcommand{\Paravalidate}{\ensuremath{\Para_{\mathrm{validation}}}}
\newcommand{\nparams}{{\ensuremath{n_\params}}}

\newcommand{\ndof}{\ensuremath{n}}
\newcommand{\Figure}[1]{Figure #1}

\newcommand{\Section}[1]{Section #1}

\newcommand{\methodname}{ROMES}%
\newcommand{\rbmatlab}{RBMatlab}
\newcommand{\methodexplanation}{ROM Error Surrogates}%
\newcommand{\stochastic}[1]{\widetilde{#1}}
\newcommand{\machine}{\quantitysurr}
\newcommand{\machiney}{\tilde y}
\newcommand{\ciVar}{{\omega}}
\newcommand{\cobs}[1]{{\ciVar_\text{validation}\left(#1\right)}}
\newcommand{\errorSurrUni}{\stochastic{\error_\text{uni}}}
\newcommand{\functionErrorNo}{d}
\newcommand{\bepsilon}{\bm{\epsilon}}
\newcommand{\functionError}[1]{\functionErrorNo(#1)}
\newcommand{\errorsurr}{\tilde \functionErrorNo}
\newcommand{\rbphi}{\psi}
\newcommand{\bphi}{\bm \phi}

\newcommand{\rvmvar}{k}
\newcommand{\rvmnum}{K}

\newcommand{\Mgreedy}{\cM_{\text{greedy}}}
\newcommand{\identityFunction}[1]{\mathrm{id}_{#1}}

\newcommand{\nrb}{p}
\newcommand{\nrbdual}{p_y}
\newcommand{\expected}[1]{\ensuremath{\mathrm{E}\eb{#1}}}
\newcommand{\prob}[1]{\ensuremath{\mathrm{P}[#1]}}
\newcommand{\effectivity}{\ensuremath{\eta}}
\newcommand{\effectivityLow}{\ensuremath{\eta_1}}
\newcommand{\effectivityHigh}{\ensuremath{\eta_2}}
\newcommand{\febasisvec}{\ensuremath{\varphi}}
\newcommand{\dirac}{\ensuremath{\delta_\text{Dirac}}}
\newcommand{\coefficientFunction}{\ensuremath{b}}

\usepackage{ifthen}

\newcommand{\setaxiscaption}[1]{\gdef\thenextaxiscaption{#1}}
\newcommand{\axiscaption}{\expandafter\thenextaxiscaption}
\newcommand{\imagefile}[1]{/home/mdrohma/latex/rb_paper/images/#1}

\newcommand{\inputaxis}[2]{{\pgfplotsset{every axis/.append style={#2}}\input{\imagefile{#1}}}}

\newcommand{\rbmedium}{18}
\newcommand{\rblarge}{62}
\newcommand{\vardefault}{1}
\newcommand{\dirdefault}{res2energy}

\ifthenelse{\equal{\withPgf}{yes}}
{}
{}

\setaxiscaption{$\backslash$axiscaption does not work if externalized!!}

\usepackage{./definitions/tikzexternal}
\ifthenelse{\equal{\atype}{siam}}{
\tikzsetexternalprefix{siam_extfigures/}}
{
\tikzsetexternalprefix{extfigures/}}

\ifthenelse{\equal{\atype}{siam}}{
\externaldocument{supplements}
}{}

\usepackage{fancyhdr}
\title{The \methodname{} method for statistical modeling of
reduced-order-model error}
\ifthenelse{\equal{\atype}{siam}}
{
\author{
Martin Drohmann\footnotemark[1]
\and
Kevin Carlberg\footnotemark[1]}
}{
}
\ifthenelse{\equal{\atype}{preprint}}
{
\author{Martin Drohmann \and Kevin Carlberg}
\address[Martin Drohmann]{Sandia National Laboratories, 7100 East Ave, Livermore, California 94550}
\email[M.~Drohmann]{mdrohma@@sandia.gov}
\address[Kevin Carlberg]{Sandia National Laboratories, 7100 East Ave, Livermore, California 94550}
\email[K.~Carlberg]{ktcarlb@@sandia.gov}
}{
}

\begin{document}

\maketitle
\begin{abstract}
  This work presents a technique for statistically modeling errors introduced
  by reduced-order models. The method employs Gaussian-process regression to
  construct a mapping from a small number of computationally inexpensive
  `error indicators' to a distribution over the true error. The variance of
  this distribution can be interpreted as the (epistemic) uncertainty
  introduced by the reduced-order model.
  To model normed errors, the method employs existing rigorous error
  bounds and residual norms as indicators; numerical experiments show that the method leads to a
  near-optimal expected effectivity
  in contrast to typical error bounds. To model errors in general outputs, the
  method uses dual-weighted residuals---which are amenable to uncertainty
  control---as indicators. Experiments illustrate that correcting the
  reduced-order-model output with this surrogate can improve prediction accuracy
  by an order of magnitude; this contrasts with existing `multifidelity
  correction' approaches, which often fail for reduced-order models and suffer
  from the curse of dimensionality.  The proposed error surrogates also lead
  to a notion of `probabilistic rigor', i.e., the surrogate bounds the error
  with specified probability.
\end{abstract}

\ifthenelse{\equal{\atype}{siam}}{
\newcommand{\slugmaster}{%
\slugger{juq}{xxxx}{xx}{x}{x--x}}
}{}

\pagestyle{fancy}
\rfoot{\vskip0.3em Preprint submitted to JUQ}

\renewcommand{\thefootnote}{\fnsymbol{footnote}}
\ifthenelse{\equal{\atype}{siam}}{
\footnotetext[1]{Sandia National
Laboratories\\7011 East Ave, MS 9159, Livermore, CA 94550, ({\tt
\{mdrohma, ktcarlb\}@sandia.gov}).  Questions, comments, or
corrections to this document may be directed to these email
addresses.}
}{}
\renewcommand{\thefootnote}{\arabic{footnote}}

\ifthenelse{\equal{\atype}{siam}}{
\begin{keywords}
model reduction, uncertainty quantification, a posteriori error estimation,
Gaussian processes, supervised machine learning
\end{keywords}
}{}

\ifthenelse{\equal{\atype}{siam}}{
\begin{AMS}
  65G99,%
  65Y20,%
  62M86%
\end{AMS}

\pagestyle{myheadings}
}{}

\ifthenelse{\equal{\atype}{preprint}}{
\subjclass{Primary 65G99, Secondary 65Y20, 62M86}}{}

\thispagestyle{plain}
\markboth{\methodname{} METHOD}{M.~DROHMANN, K.~CARLBERG}

\section{Introduction} \label{sec:introduction} As computing power increases,
computational models of engineered systems are being employed to answer
increasingly complex questions that guide decision making, often in
time-critical scenarios. It is becoming essential to rigorously quantify and account
for both aleatory and epistemic uncertainties in
these analyses. 
Typically, the high-fidelity computational model can
be viewed
as providing a (costly-to-evaluate) mapping between system \emph{inputs}
(e.g., uncertain parameters, decision variables) and system \emph{outputs}
(e.g., outcomes, measurable quantities). For example, data
assimilation employs collected sensor data (outputs) to update the
distribution of uncertain
parameters (inputs) of the model; doing so via
Bayesian inference requires sampling from the posterior distribution, which
can entail thousands of forward model simulations. The computational 
resources (e.g., weeks on a supercomputer) required for 
large-scale simulations preclude such high-fidelity models from being feasibly
deployed in such scenarios. 

To avoid this bottleneck, analysts have turned to surrogate
models that approximate the input--output map of the high-fidelity model, yet
incur a fraction of their computational cost. However, to be rigorously
incorporated in uncertainty-quantification (UQ) contexts, it is critical to quantify the
\emph{additional uncertainty} introduced by such an approximation. For example, 
Bayesian inference aims to sample from the posterior distribution
\begin{equation}
\prob{\inputs|\outmeas}
\propto\prob{\inputs}\prob{\outmeas|\inputs},
\end{equation}
where $\bmu \in \Para\subset\RR^\nparams$ denote system inputs, $\outmeas\in\RR$ denotes
the measured output,\footnote{This work considers one output for
notational simplicity. All concepts can be straightforwardly extended to multiple
outputs. The numerical experiments treat the case of multiple outputs.}
$\prob{\inputs}$ represents the prior, and $\prob{\outmeas|\inputs}$ denotes the likelihood
function. Typically, the measured output is modeled as
$\outmeas = \out(\inputs) + \varepsilon$, where $\out:\Para\rightarrow \RR^{}$
denotes the
outputs predicted by the high-fidelity model for inputs $\inputs$,
and
$\varepsilon$ is a random variable representing measurement noise.
Sampling from this posterior distribution (e.g., via Markov-chain Monte--Carlo
or importance sampling) is costly, as each sample requires at least
one evaluation of the high-fidelity input--output map $\inputs\mapsto\out$
that appears in the likelihood function.

When a surrogate model is employed, the measured output becomes $\outmeas =
\outsurr(\inputs) + \outbias(\inputs) + \varepsilon$, where
$\outsurr:\Para\rightarrow \RR$ denotes the output predicted by the surrogate
model, and $\outbias:\Para\rightarrow \RR^{}$ represents the
\emph{surrogate-model output error} or bias. In this case, posterior sampling 
requires only evaluations of the \emph{surrogate-model} input--output map
$\inputs\mapsto\outsurr $---which is computationally inexpensive---as well as
evaluation of the surrogate-model error $\outbias(\inputs)$, which is not
precisely known in practice. As such, it can be considered a source of
\emph{epistemic uncertainty}, as it can be reduced in principle by employing the
original high-fidelity model (or a higher fidelity surrogate model).
The goal of this work is to construct a
statistical model of this surrogate-model error
$\stochastic{\outbias}(\inputs)$ that is
1) cheaply computable, 2) exhibits low variance (i.e., introduces minimal
epistemic uncertainty), and 3) whose distribution can be numerically
validated.

Various approaches have been developed for different surrogate models to
quantify the surrogate error $\outbias(\inputs)$.  Surrogate models can be
placed into three categories \cite{Eldred2004}: 1) data fits, 2) lower-fidelity models, and 3) reduced-order models. Data
fits employ supervised machine-learning methods (e.g.,
Gaussian processes, polynomial interpolation \cite{FSK08}) to directly model the
high-fidelity input--output map. 
Within this class of surrogates, it is
possible to statistically model the error for stochastic-process data fits, as 
a prediction for inputs $\inputs$ yields a mean $\outsurr(\inputs)$ and a mean-zero distribution
$\outbias(\inputs)$ that can be associated with epistemic uncertainty.
While such models are (unbeatably) fast to query and non-intrusive to
implement,\footnote{Their construction requires only black-box evaluations of
the input--output map of the high-fidelity model.} they suffer from the curse
of dimensionality and lack access to the underlying model's physics, which
can hinder predictive robustness.  

Lower-fidelity models simply replace the high-fidelity model with a
`coarsened' model obtained by neglecting physics, coarsening the mesh, or
employing lower-order finite elements, for example.  While such models remain
physics based, they often realize only modest computational
savings. For such problems, `multifidelity correction' methods have
been developed, primarily in the optimization context. These techniques
model the mapping $\inputs\mapsto\outbias$ using a data-fit surrogate;
they either enforce `global' zeroth-order consistency between the corrected
surrogate prediction and the high-fidelity prediction at training points
\cite{gano2005hybrid,huang2006sequential,march2012provably,rajnarayan2008multifidelity,NE12}, or `local' first- or
second-order consistency at trust-region centers \cite{Alexandrov2001,
Eldred2004}.  Such approaches tend to work well when the surrogate-model error
exhibits a lower variance than the high-fidelity response \cite{NE12} and the input-space dimension is small.

Reduced-order models (ROMs) employ a projection process to reduce the
state-space dimensionality of the high-fidelity computational model. Although
intrusive to implement, such physics-based surrogates often lead to more
significant
computational gains than lower-fidelity models, and higher robustness than
data fits. For such models, error analysis has been limited primarily to
computing rigorous \textit{a posteriori} error bounds
$\outbound(\bmu)$ satisfying $\outerror{(\bmu)}\leq\outbound(\bmu)$
\cite{BM12, GP05, RHP07}.  Especially for nonlinear problems, however, these
error bounds are often highly ineffective, i.e., they overestimate the actual error
by orders of magnitude \cite{DHO12}.  To overcome this shortcoming and obtain
tighter bounds, the ROM must be equipped with complex machinery that both increases the computational burden \cite{WSH12,
HRSP06} and is intrusive to implement (e.g., reformulate the discretization of the high-fidelity model
\cite{UP13,YP13}).  Further, \emph{rigorous} bounds are not directly useful
for uncertainty quantification (UQ) problems, where a statistical error model that is
unbiased, has low variance, and is stochastic is more useful.  Recent work
\cite[Section IV.D]{NE12} has
applied multifidelity correction to ROMs. However, the method did not succeed
because the ROM error is often
a highly oscillatory function of the inputs and therefore typically exhibits a
\emph{higher} variance than the high-fidelity response.

  \begin{table}[t]
    \centering
      \tikzsetfigurename{feature_matrix}
      \begin{tikzpicture}
        \newcommand{\bad}{}
        \newcommand{\good}{}
        \newcommand{\medium}{}
        \node {
        \begin{tabular}{lcccc}
          \toprule
           & ROM & data& multifidelity & ROM + \\
           & & fits & correction& \methodname \\
          \midrule
          non--intrusive & \bad $\times$ & \good $\checkmark$ & \good
          $\checkmark$ & \bad $\times$ \\
          output-error correction & \bad $\times$ & N/A & \good $\checkmark$ &  \good $\checkmark$ \\
          rigorous error bounds & \good $\checkmark$ & \bad $\times$ & \bad $\times$ & \medium $(\checkmark)^*$ \\
          tight error bounds & \medium $(\checkmark)^{\dagger}$ & \bad $\times$ &\bad $\times$ & \good $\checkmark$ \\
          \bottomrule
        \end{tabular}
        };
      \end{tikzpicture}
      \begin{minipage}{0.7\textwidth}
        \footnotesize
        \begin{\itemname}
          \item[*] probabilistically rigorous\\[-1.7em]
          \item[$\dagger$] good effectivity can only be obtained with very
					intrusive methods.
        \end{\itemname}
      \end{minipage}
    \caption{Features of different surrogate models \label{tab:featureMatrix}
     }
  \end{table}
In this paper, we introduce the \methodexplanation\ (\methodname{}) method that
aims to combine the utility of multifidelity correction with the computational
efficiency and robustness of reduced-order modeling.  Table \ref{tab:featureMatrix} compares
the proposed approach with existing surrogate-modeling techniques. Similar to the multifidelity-correction
approach, we aim to model the ROM error $\outbias$ using a data-fit surrogate.
However, as directly
approximating the mapping $\inputs\mapsto\outbias$ is ineffective for ROMs,
we instead exploit the following key observation: ROMs often generate a small
number of physics-based,
cheaply computable \emph{error indicators} $\stateindicator:\Para\rightarrow \RR^{\nindicator}$ that correlate strongly with the true error
$\outbias(\bmu)$. Examples of indicators include the
residual norm, dual-weighted residuals,
and the rigorous error bounds discussed above. To this end, \methodname\
approximates
the \emph{low-dimensional, well-behaved} mapping $\stateindicator(\bmu)\mapsto
\outbias(\bmu)$ using Gaussian-process regression, which is a
stochastic-process data-fit method. Note that \methodname\ constitutes a
generalization of the multifidelity correction approach, as the inputs (or
features) of the
error model can be any user-defined
error indicator---they need not be the system inputs $\bmu$.
\Figure \ref{fig:OURSchema} depicts the propagation of information for the
proposed method.

In addition to constructing an error surrogate for the system outputs,
\methodname{} can also be used to construct a statistical model for the norm
of the error in the system state. Further, \methodname{}
can be used to generate error bounds with `probabilistic rigor', i.e.,  an error
bound that overestimates the error with a specified probability.

  \begin{figure}[htp]
    \centering
      \tikzsetfigurename{ROMSESSchema}
      \begin{tikzpicture}
        {
        \tikzset{every node/.append style={draw=black}}
        \tikzset{split/.style={rectangle split, rectangle split parts=2}}
        \tikzset{arrow/.style={thick, ->}}
        \begin{scope}
          \node (input) [split] {\em Input \nodepart{two} $\bmu$};
          \node (system) at ($(input.east)+(2.8cm,0)$) [anchor=west, text width=4.5cm, align=center]
          { {\em Reduced-order model}
          };
          \node (output) at ($(system.east)+(2.8cm,0)$) [anchor=west, split]
          {\em Output \nodepart{two} $\redout(\bmu)\defeq \out(\redstate(\bmu))$};
          \node (indicator) at ($(input.south west)-(0, 0.5cm)$) [anchor=north west, split] {\em Error indicators \nodepart{two} $\stateindicator(\bmu)$};

          \node (romessystem) at ($(system.south)-(0,1.1cm)$) [anchor=north] { {\em \methodname}};

          \node (bound) at ($(output.south east)-(0, 0.4cm)$) [anchor=north east, split]
          {\em Error surrogates \nodepart{two} error
					$\stochastic{\stateerror{}}(\indicators(\bmu))$ or output error $\stochastic{\outbias}(\indicators(\bmu))$};

          \draw[arrow] (input) -- (system);
          \draw[arrow] (system) -- (output);
          \draw[arrow] (system) -- (indicator);
          \draw[arrow] (indicator) -- (romessystem);
          \draw[arrow] (romessystem) -- (bound);
        \end{scope}
        }
      \end{tikzpicture}
      \caption{\methodname{} method. The output quantities of interest can be `corrected'
      by adding the ROM error surrogate to the ROM
      output prediction, i.e., $\out(\bmu)\approx
      \outcorrected(\bmu)\defeq \redout(\bmu)+\stochastic{\outbias}(\indicators(\bmu))$.
      }\label{fig:OURSchema}
  \end{figure}

Next, \Section{\ref{sec:probFormulation}} introduces the problem formulation
and provides a general but brief introduction to model reduction.
  In \Section{\ref{sec:notation}}, we introduce the \methodname{}
  method, including its objectives, ingredients, and some choices of these
	ingredients for particular errors.
  \Section{\ref{sec:gpbasics}} briefly
  summarizes the Gaussian-process kernel method \cite{RW06} and the relevance vector
  machine \cite{Tipping2001}, which are the two machine-learning algorithms we
	employ to construct the \methodname{} surrogates.
  However, the \methodname{}
  methodology does not rely on these two techniques, as any 
  supervised machine learning algorithm that generates a stochastic process
	can be used, as long as it generates a statistical model that meets
	the important conditions described in \Section{\ref{sec:statModel}}.
  \Section{\ref{sec:experiments}} analyses
  the performance of the method when applied with the reduced-basis method to
	solve Poisson's equation in two dimensions using nine system inputs. The
	method is also compared with existing rigorous error bounds for normed
	errors, and with the multifidelity correction approach for errors in general
	system outputs.

  For additional information on the reduced-basis method, including the
  algorithms to generate the reduced-basis spaces and the computation of error bounds,
  we refer to the supplementary \Section{\ref{sec:rb}}.

\section{Problem formulation}\label{sec:probFormulation}
This section details aspects of the high-fidelity and reduced-order models
that are important for the \methodname{}
surrogates.  We begin with a formulation of the high-fidelity model in
Section \ref{sec:problem:full} and the reduced-order model in
Section \ref{sec:problem:reduced}. Finally, we elaborate on the errors
introduced by the model-reduction process and possible problems with their
quantification in Section \ref{sec:problem:errors}.

\subsection{High-fidelity model\label{sec:problem:full}}
Consider solving a parameterized systems of equations
\begin{equation}\label{eq:FOMdyn}
\residual(\state;\params) = 0,
\end{equation}
where
   $\state:\Para\rightarrow\RR^{\ndof}$ denotes the state implicitly defined by Eq.\ \eqref{eq:FOMdyn},
  $\params\in\paramDomain\subset\RR^{\nparams}$ denotes the system inputs, and $\residual:\RR^{\ndof}\times\RR^{\nparams}\rightarrow \RR^{\ndof}$ denotes the
  residual operator. This model is appropriate for stationary problems, e.g.,
those arising from the finite-element discretization of elliptic PDEs.
For simplicity, assume we are interested in computing a single
output
\begin{equation} \label{eq:outdyn}
\outputs(\params)  \defeq \outputfun(\state(\params))
\end{equation}
with $\outputs:\RR^{\nparams}\rightarrow\RR$ and
$\outputfun:\RR^{\ndof}\rightarrow\RR^{}$.

When the dimension $\ndof$ of the high-fidelity model is `large', computing
the system outputs $\outputs$ by first
solving Eq.~\eqref{eq:FOMdyn} and subsequently applying
Eq.~\eqref{eq:outdyn} can be prohibitively expensive. This is particularly
true for many-query problems arising in UQ such as Bayesian inference, which
may require
thousands of input--output map evaluations $\params\mapsto\outputs$.

\subsection{Reduced-order model\label{sec:problem:reduced}}
Model-reduction techniques aim to reduce the burden of solving Eq.\
\eqref{eq:FOMdyn} by employing a projection process.
First, they execute a computationally expensive \emph{offline stage}
(e.g., solving Eq.~\eqref{eq:FOMdyn} for a training
set $\params\in\Paratrain\subset\paramDomain$) to construct 1) a
low-dimensional trial basis (in matrix form) $\rightbasis\in
\RR^{\ndof\times\nrb}$ with $\nrb\ll\ndof$ that (hopefully) captures the behavior of
the state $\state$ throughout the parameter domain $\paramDomain$, and 2) an associated test basis $\leftbasis\in\RR^{\ndof\times\nrb}$.
Then, during the computationally inexpensive \emph{online stage}, these methods
approximately solve Eq.~\eqref{eq:FOMdyn} for arbitrary $\params\in
\paramDomain$ by searching for solutions
in the trial subspace $ \range{\rightbasis}\subset\RR^{\ndof}$ and
enforcing orthogonality of the residual $\residual$ to the test subspace
$\range{\leftbasis}\subset\RR^{\ndof}$:
\begin{equation}\label{eq:ROMdyn}
\leftbasis^t\residual(\rightbasis\redstatecoords;\params) = 0.
\end{equation}
Here, the state is approximated as
$\redstate(\params)\defeq\rightbasis\redstatecoords(\params)$ and the reduced state
$\redstatecoords(\params)\in\RR^{\nrb}$ is implicitly defined by
Eq.~\eqref{eq:ROMdyn}.
The ROM-predicted output is then
$\redout(\bmu)\defeq\outputfun\left(\redstate(\params);\params\right)$.

When the residual operator is nonlinear in the state or non-affine in the
inputs, additional
complexity-reduction approximations such as empirical interpolation
\cite{barrault2004eim,grepl2007erb}, collocation
\cite{LeGresleyThesis,astrid2007mpe,ryckelynck2005phm}, discrete
empirical interpolation
\cite{chaturantabut2010journal,galbally2009non,drohmannEOI},  or gappy proper orthogonal
decomposition (POD) \cite{carlbergGappy,carlbergJCP} are required to ensure that computing the low-dimensional
residual $\leftbasis^t\residual$ incurs an $\ndof$-independent
operation count. In this case, the residual is approximated as $\residualhyper\approx\residual$ and the reduced-order equations become
\begin{equation}\label{eq:ROMdynhyper}
\leftbasis^t\residualhyper(\rightbasis\redstatecoords;\params) = 0.
\end{equation}
When the output operator is nonlinear and the vector $\partial g/\partial
\state$ is dense, approximations in the output
calculation are also required to ensure an $\ndof$-independent operation count.

Section \ref{sec:rb} describes in detail the construction of
a reduced-order model using the reduced-basis method applied to a
parametrically coercive, affine, linear, elliptic PDE.

\subsection{Reduced-order-model error bounds}\label{sec:problem:errors}
One is typically interested in quantifying two types of error incurred by model
reduction: the state-space error $\stateerrorvec{}(\params)\defeq
\state(\params)- \redstate(\params) \in\RR^{\ndof}$ and the output error
    $\outbias{(\bmu)} \defeq   \out(\bmu)-\redout(\bmu) \in\RR{}$. In
    particular, many ROMs are equipped with computable, rigorous error bounds
    for these quantities \cite{PR07, GMNP06, NR09, CTU09, DHO12}:
\begin{equation}\label{eq:errorBounds}
\statebound{(\bmu)}\geq\stateerror{(\bmu)},\qquad
\outbound{(\bmu)}\geq\outerror{(\bmu)}
\end{equation}
In cases, where the norm of the residual operator can be estimated tightly, lower
bounds also exist:
\begin{equation}\label{eq:lowerErrorBounds}
\statelowerbound{(\bmu)}\leq\stateerror{(\bmu)},\qquad
\outlowerbound{(\bmu)}\leq\outerror{(\bmu)}.
\end{equation}
The performance of an upper bound is usually quantified by its {\em
effectivity}, i.e., the factor by which it overestimates the true error
\begin{equation}
  \eta^\state(\bmu) \defeq  \frac{\statebound{(\bmu)}}{\stateerror{(\bmu)}}\geq 1,
  \qquad
  \eta^\out(\bmu) \defeq  \frac{\outbound{(\bmu)}}{\outerror{(\bmu)}}\geq 1.
  \label{eq:effectivityBoundDefinition}
\end{equation}
The closer these values are to $1$, the `tighter' the bound. For coercive PDEs,
these effectivities can be controlled by choosing a tight lower
bound of the coercivity constant.

While this can be easily
accomplished for stationary, linear problems, it is difficult to
find tight lower bounds in almost all other cases. In fact, the resulting
bounds
  often overestimate the error by orders of magnitude \cite{RHP07, DHO12}. Because effectivity is 
critically important in practice, various efforts
have been undertaken to improve the tightness of the bounds. Huynh et al.~\cite{HRSP06}
developed the successive constraint method for this purpose; the
method approximates the
coercivity-constant lower bounds by solving small linear programs online,
which depend on
additional expensive offline computations. Alternatively, Refs.~\cite{UP13, YP13} reformulate the
entire discretization of time-dependent problems using a space--time
method that improves the error bounds by incorporating solutions to dual problems. 
Another approach \cite{WSH12} aims to
approximate the coercivity constant by eigenvalue analysis of the reduced
system matrices.  These methods all bloat the offline and the online
computation time and often incur intrusive changes to the high-fidelity-model
implementation.

 Regardless of effectivity, rigorous bounds satisfying inequalities
\eqref{eq:errorBounds} are not directly useful for quantifying the epistemic
uncertainty incurred by employing the ROM. Rather, a \emph{statistical model}
that reflects our knowledge of these errors would be more appropriate. For
such a model, the
mean of the distribution would provide an expected error; the variance would
provide a notion of epistemic uncertainty.
The most straightforward way to achieve this
would be to model the error as a uniform probability distribution on
an interval whose boundaries correspond to
the lower and upper bounds. Unfortunately, such an approach leads to wide,
uninformative intervals when the bounds suffer from poor effectivity; this will be demonstrated in the numerical experiments.

Instead, we exploit the following observation: error bounds tend to be
strongly correlated with the true error.
\Figure{\ref{fig:correlations}} depicts this observed
structure for a reduced-basis ROM
  applied to an elliptic PDE (see \Section \ref{sec:experiments:thermalblock}
	for details). On a logarithmic
scale, the true error exhibits a roughly linear dependence on both the bound and the
residual norm, and the variance of the data is fairly constant. As
will be shown in \Section \ref{sec:experiments:comparison:Eldred}, employing a multifidelity correction
approach wherein the error is modeled as a function of the inputs
$\params$ does not work well for this example, both because the input-space
dimension is large (nine) and the error is a highly oscillatory function of these
inputs.
  \begin{figure}
    \begin{center}
      \tikzsetfigurename{correlations}
      \begin{tikzpicture}
          \pgfplotsset{
            tick label style={font=\small},
            label style={font=\small, text width=4.0cm, align=center},
            legend style={font=\small},
            width=6cm,
          }

          \inputaxis{\dirdefault/correlation_1d_bs_\rbmedium_v1_tikz.tex}{%
              legend to name=leg:correlation1d,
              legend entries={{residuals}, {error bounds}},
              name=axis_correlation1d}

          \inputaxis{\dirdefault/correlation_2d_bs_\rbmedium_v1_tikz.tex}{%
            at={($(axis_correlation1d.east)+(3cm,0)$)}, anchor=west}

          \node at ($(axis_correlation1d.south east)-(0.1cm,-0.1cm)$) [anchor=south east]
            {\pgfplotslegendfromname{leg:correlation1d}};

      \end{tikzpicture}
      \caption{Relationship between RB error bounds $\statebound$, residual
      norms $\|\residual(\rightbasis\redstatecoords;\params)\|$,
      and the true state-space errors $\enstateerror{}$, visualized by evaluation of
      $200$ random sample points in the input space. Here, 
      $\enorm{\cdot}$ denotes the energy norm defined in
      Section~\ref{sec:experiments:thermalblock}.
      }\label{fig:correlations}
    \end{center}
  \end{figure}

Therefore, we propose constructing a stochastic process that maps such
\emph{error indicators} to a \emph{random variable for the error}. For
this purpose, we employ Gaussian-process regression. The approach leverages one strength of ROMs
compared to other surrogate models: ROMs generate strong `physics-based' error indicators (e.g., error bounds) in addition
to output predictions. The next section describes the proposed method.

\section{The \methodname{} method} \label{sec:notation}
    The objective of the \methodname{} method is to construct a statistical
    model of the deterministic, but generally unknown ROM error
    $\error(\params)$ with $\error:\Para \to \RR$ denoting an
    $\RR{}$-valued error that may represent the norm of the state-space error
    $\stateerror{}$, the output error $\outbias$, or its absolute value
    $|\outbias|$, for example.  The distribution of the random variable
    representing the error should reflect our (epistemic) uncertainty about
    its value. We assume that we can employ a set of training points
    $\error(\paramTrain{n})$, $n=1,\ldots,\nTrain$ to construct this model.

  \subsection{Statistical model}\label{sec:statModel}

  Define a probability space $(\Omega,\mathcal F,P)$. We seek to approximate
  the deterministic mapping $\quantity:\params\mapsto
  \functionError{\error(\params)}$ by a stochastic mapping 
  $\quantitysurr:\indicators(\params)\mapsto\errorsurr$ with
  $\errorsurr:\Omega\rightarrow\RR{}$ a real-valued random variable. Here,
  $\functionErrorNo:\RR{}\rightarrow\RR{}$ is an invertible
  \emph{transformation function} (e.g., logarithm) that
can be specified to facilitate construction of the statistical model. We
can then interpret the statistical model of the error as a random variable
$\tilde \delta :\indicators\mapsto
\functionErrorNo^{-1}(\quantitysurr(\indicators))$.

Three ingredients must be selected to construct this mapping $\quantitysurr$:
1) the error indicators $\indicators$, 2) the transformation function
$\functionErrorNo$, and 3) the
methodology for constructing the statistical model from the training data. We
will make these choices such that the stochastic mapping satisfies the
following conditions:
\begin{enumerate} 
  \item The indicators $\indicators(\params)$ are \emph{cheaply computable}
	and \emph{low dimensional}
  given any $\params\in\Para$. In practice, they should also incur a
  reasonably small implementation effort,
  e.g., not require modifying the underlying high-fidelity model.
  \item\label{cond:lowvar} The mapping $\quantitysurr$
  exhibits \emph{low variance}, i.e., 
  $\expected{\left(\quantitysurr(\indicators(\params)) - 
  \expected{\quantitysurr(\indicators(\params))}\right)^2}$ is `small' for all
  $\params\in\Para$. This ensures that little additional epistemic
uncertainty is introduced.
  \item\label{cond:validated} The mapping $\quantitysurr$ is \emph{validated}:
 \begin{equation} \label{eq:ciValidate}
 \cobs{\ciVar}\approx \ciVar,\quad \forall \ciVar\in\left[0,1\right),
  \end{equation} 
where $\cobs{\ciVar}$ is the frequency with which validation data lie in the
$\ciVar$-confidence interval predicted by the statistical model
\begin{equation}\label{eq:validatedModelA}
\cobs{\ciVar}\defeq\frac{\card{\{\bmu\in\Paravalidate\ |\ 
\functionError{\error(\params)}\in\confidenceInterval{\ciVar}{\params}\}}}{\card{\Paravalidate}}.
\end{equation}

Here, the validation set $\Paravalidate\subset\Para$ should not include any of
the
points 
$\paramTrain{n}$, $n=1,\ldots,\nTrain$
employed to train the error surrogate, and the confidence interval $\confidenceInterval{\ciVar}{\params}
\subset \RR{}$, which is centered at the mean of
$\quantitysurr(\indicators(\params))$,
 is defined for all
$\params\in\Para$ such that 
\begin{equation}
  \prob{\quantitysurr(\indicators(\params)) \in
	{\confidenceInterval{\ciVar}{\params)}}} = \ciVar.
\end{equation}
	In essence, validation assesses whether or not the data do indeed
	behave as random variables with probability distributions predicted by the
	statistical model.
\end{enumerate}
The next section describes the proposed methodology for selecting indicators
$\indicators$ and transformation function $\functionErrorNo$. For constructing the mapping
$\machine$ from training points, we will employ the two {supervised
machine learning algorithms} described in \Section \ref{sec:gpbasics}: the
Gaussian process (GP) kernel method and the relevance vector machine (RVM).
  Note that these are merely guidelines for model construction,
  as there are usually no strong analytical tools to prove
  that the mapping behaves according to a certain probability distribution.
  Therefore, any choice must be computationally validated according to
	condition \ref{cond:validated} above.

\subsection{Choosing indicators 
and transformation function}
\label{sec:romes:stochasticMapping}

The class of multifidelity-correction algorithms can be cast within the
framework proposed in Section \ref{sec:statModel}.  In particular, when a
stochastic process is used to model additive error, these methods are
equivalent to the proposed construction with ingredients $\delta =
\outbias$, $\indicators = \params$, and $\functionErrorNo =
\identityFunction{\RR{}}$ with $\identityFunction{\RR{}}(x) = x$, $\forall
x\in\RR{}$ the identity function over $\RR{}$.  However, as previously
discussed, the mapping $\bmu \mapsto \outbias$ can be highly oscillatory and
non-smooth for reduced-order models; further, this approach is infeasible for
high-dimensional input spaces, i.e., $\nparams$ large. This was shown  by Ng
and Eldred \cite[Section IV.D]{NE12}; we also demonstrate this in the
numerical experiments of Section \ref{sec:experiments:comparison:Eldred}.

Note that all indicators and errors proposed in this section should be scaled (e.g., via
linear transformations) in practice such that they exhibit roughly the same
range.  This `feature scaling' task is common in machine-learning and is
particularly important when the \methodname{} surrogate employs multiple
indicators.

\subsubsection{Normed and compliant outputs}\label{sec:errorboundindicators}

As discussed in Section \ref{sec:problem:errors}, many ROMs are equipped with
bounds for normed errors.  Further, there is often a strong, well-behaved
relationship between such error bounds and the normed error (see Figure
\ref{fig:correlations}).  In the case of
compliant outputs, the error is always non-negative, i.e., $\outbias =
|\outbias| $ (see Section \ref{sec:ellipticProblem:errorBounds}), so we can treat this error as a normed error.

To this end, we propose employing error bounds as indicators 
for the errors in the compliant output $|\outbias|$ and in the state
$\stateerror{}$. However, because the bound effectivity often lies in a small
range (even for a large
range of errors) \cite{RHP07}, employing a
logarithmic transformation is appropriate. To see this, consider a case where
the effectivity $\effectivity$ of the error bound, defined as 
\begin{equation}
  \effectivity(\bmu) \defeq  \frac{\bound{(\bmu)}}{\error{(\bmu)}}\geq 1,\quad
  \forall\bmu\in\Para,
  \label{eq:effectivityBoundDefinitionGen}
\end{equation}
lies within a range
$\effectivityLow\leq\effectivity(\bmu)\leq\effectivityHigh$, $\forall
\bmu\in\Para$. Then the relationship between the error bound and the true
error is
\begin{gather} 
\frac{\bound{(\bmu)}}{\effectivityLow}\geq\error{(\bmu)}\geq\frac{\bound{(\bmu)}}{\effectivityHigh}\\
\log{\bound{(\bmu)}}-\log{\effectivityLow}\geq\log{\error{(\bmu)}}\geq\log{\bound{(\bmu)}}-\log{\effectivityHigh}
\end{gather} 
for all $\bmu\in\Para$.  In this case, one would expect an affine model
  mapping $\log{\bound{(\bmu)}}$ to $\log{\error{(\bmu)}}$ with constant
  Gaussian noise to accurately capture the relationship. So, employing
  a logarithmic transformation permits the use of simpler surrogates that
  assume a constant variance in the error variable.
Therefore, we propose employing $\indicators = \log{\bound}$
and $\functionErrorNo = \log$ for statistical modeling of normed errors and
compliant outputs.

A less computationally expensive candidate for an indicator is simply the
logarithm of the residual norm
$\indicators = \log{\residualnorm}$, where $\residualnorm$ is the Euclidean
norm of the residual vector
  \begin{equation}
    \residualnorm(\bmu) \defeq
		\norm{\residual(\rightbasis\redstatecoords(\params);\params)}_2.
    \label{eq:residualExampleDef}
  \end{equation}
For more information on the efficient computation of
\eqref{eq:residualExampleDef}, we refer to Section \ref{sec:offlineOnlineResidual}.
One would expect this choice of indicator to produce a similar model to that
produced by the logarithm of the error bound: the error bound is often equal
to the residual norm divided by a (costly-to-compute) coercivity-constant
lower bound
(see Section \ref{sec:ellipticProblem:errorBounds}). Further, employing the
residual norm leads to a model that is less sensitive to strong variations in
this approximated lower bound.

 Returning to the example depicted in \Figure{\ref{fig:correlations}}, the
 relationship between the error bound and the energy norm of the state error
 in log-log space
 is roughly linear, and the variance is relatively small. The same is true of the
 relationship between the (computationally cheaper) residual norm  and the
 true error. As expected, these relationships can be accurately modeled as a
 stochastic process with a linear mean function and constant variance (more
 details in \Section \ref{sec:experiments:thermalblock}).  Here, strong
 candidates for \methodname{} error indicators include $\indicators_1(\bmu)\defeq \residualnorm(\bmu),
 \indicators_2(\bmu)\defeq \bound{(\bmu)}$, and $\indicators_3(\bmu)\defeq \rb{\residualnorm(\bmu),
 \bound{(\bmu)}}$. In the experiments in Section \ref{sec:experiments},
 we will consider the first choice, which is the least expensive and 
 intrusive option, yet leads to excellent results.  For cases where the data are
 less well behaved, more error indicators can be included,
  e.g., linear combinations of inputs or the output prediction.

Unfortunately, this set of \methodname{} ingredients is not applicable to
errors in \emph{general} outputs of interest
because the logarithmic transformation function assumes strictly positive errors.
The next section presents a strategy for handling this general case.

  \begin{remark}[Log-Normal distribution]
    \label{rem:lognormal}
    In the case where $\functionErrorNo=\log$, the error models
    $\stochastic{\error}(\bmu)$, $\bmu\in\Para$ are random
    variables with log-normal distribution. If
    one is interested in the most probable
    error, one might think to use the expected value of
    $\stochastic{\error{}}$. However,
    the maximum of the probability distribution function of a log-normally
    distributed random variable is defined by its {\em mode}, which is less
    than the expected value. We therefore use
    $
    \mathrm{mode}(\stochastic{\error{}})
    $
    if scalar values for the estimation of the output error or the reduced state error
    are required.
  \end{remark}
  \subsubsection{General outputs}\label{sec:dualWeightedResiduals}

This section describes the \methodname{} ingredients we propose for modeling the
error $\outbias$ in a general output 
$\outputs(\params)  \defeq \outputfun(\state(\params))$.
Dual-weighted-residual approaches are commonly adopted for approximating
general output errors
in the context of \textit{a posteriori}  adaptivity 
\cite{estep1995posteriori,
babuvska1984post,becker1996weighted,rannacher1999dual,venditti2000adjoint,lu2005posteriori}, model-reduction adaptivity
\cite{carlberg2014adaptive}, and model-reduction error
estimation \cite{PRVMMPT02, CTU09, UP13,meyer2003efficient}.
The latter references compute adjoint solutions in order to improve the
accuracy of ROM output-error bounds. The computation of these adjoint
solutions entails a low-dimensional linear solve; thus, they are efficiently
computable and can potentially serve as error indicators for the \methodname{}
method.

The main idea of dual-weighted-residual error estimation is to approximate the
output error to first-order using the solution to a dual problem. For
notational simplicity in this section, we drop dependence on the inputs
$\params$.

To begin, we approximate the output arising from the (unknown) high-fidelity state
$\state$ to first order about the ROM-approximated state $\rightbasis\redstatecoords$:
\begin{equation} \label{eq:outputFirst}
\outputfun\left(
\state\right)
\approx
\outputfun\left(
\rightbasis\redstatecoords\right) + \frac{\partial
\outputfun}{\partial \state}\left(
\rightbasis\redstatecoords\right)\left(\state
- \rightbasis\redstatecoords\right)
\end{equation} 
 with $\outputfun:\RR^{\ndof}\rightarrow\RR{}$ and $\frac{\partial
\outputfun}{\partial \state}:\RR^{\ndof}\rightarrow \RR^{1\times\ndof}$. Similarly, we can approximate the residual to first order about the
 approximated state as
 \begin{equation}\label{eq:resFirstEq}
 0=\rightbasisfineT\modresidual\left(
\rightbasisfine\redstatefine\right)\approx 
\rightbasisfineT\modresidual\left(
\rightbasiscoarse\redstatecoarse\right) + 
\rightbasisfineT\frac{\partial\modresidual}{\partial\state}\left(
\rightbasiscoarse\redstatecoarse\right)\rightbasisfine
\left(\redstatefine
- \coarsetofine\redstatecoarse\right),
\end{equation}
where $\residual:\RR^{\ndof}\rightarrow \RR^{\ndof}$ with
$\frac{\partial\modresidual}{\partial\state}:\RR^{\ndof}\rightarrow\RR^{\ndof\times\ndof}$.
Solving for the error yields
\begin{equation}\label{eq:resFirst}
\left(\redstatefine
- \coarsetofine\redstatecoarse\right)\approx
-
\left[\rightbasisfineT\frac{\partial\modresidual}{\partial\state}\left(
\rightbasiscoarse\redstatecoarse\right)\rightbasisfine\right]^{-1}
\rightbasisfineT\modresidual\left(
\rightbasiscoarse\redstatecoarse\right).
  \end{equation} 
Substituting \eqref{eq:resFirst} in \eqref{eq:outputFirst} leads to
\begin{equation} \label{eq:outputFinal}
\outputfun\left(
\rightbasisfine\redstatefine\right)
-
\outputfun\left(
\rightbasiscoarse\redstatecoarse\right)
\approx
 \redadjointfineT\rightbasisfineT\modresidual\left(
\rightbasiscoarse\redstatecoarse\right),
\end{equation} 
where the dual solution $\redadjointfine\in\RR^{\nrbfine}$ satisfies
\begin{equation}\label{eq:adjointFine}
\rightbasisfineT\frac{\partial\modresidual}{\partial\state}^t(
\rightbasiscoarse\redstatecoarse)\rightbasisfine\redadjointfine =
-\rightbasisfineT\frac{\partial\outputfun}{\partial\state}^t\left(
\rightbasiscoarse\redstatecoarse\right).
\end{equation}
Approximation \eqref{eq:outputFinal} is first-order accurate; therefore, it is
exact in the case of linear outputs and a linear residual operator. In the general
nonlinear case, this approximation is accurate in a neighborhood of the
ROM-approximated state $\rightbasiscoarse\redstatecoarse$.

Because we would like to avoid high-dimensional solves, we
approximate $\redadjointfine$ as the reduced-order dual solution
$\redadjointcoarsetofine\defeq\rightbasisdualcoarse\redadjointcoarse\in\RR^{\ndof}$, where
$\redadjointcoarse$ satisfies
\begin{equation}\label{eq:adjointCoarse}
\rightbasisdualcoarseT\frac{\partial\modresidual}{\partial\state}^t(
\rightbasiscoarse\redstatecoarse)\rightbasisdualcoarse\redadjointcoarse =
\rightbasisdualcoarseT\frac{\partial\outputfun}{\partial\state}^t\left(
\rightbasiscoarse\redstatecoarse\right),
\end{equation}
and $\rightbasisdualcoarse\in\RR^{\ndof\times\nrbdual}$ with
$\nrbdual\ll\ndof$ is a reduced basis (in matrix form) for the dual system. Section \ref{sec:supp:dual}
provides details on the construction of $\rightbasisdualcoarse$ for elliptic PDEs.
Substituting the approximation $\redadjointcoarsetofine$ for $\redadjointfine$
in \eqref{eq:outputFinal} yields a cheaply computable error estimate
 \begin{equation}\label{eq:outputFinalcheap} 
\outputfun\left(
\rightbasisfine\redstatefine\right)
-
\outputfun\left(
\rightbasiscoarse\redstatecoarse\right)
\approx
\redadjointcoarsetofine^t\modresidual\left(
\rightbasiscoarse\redstatecoarse\right).
  \end{equation} 
This relationship implies that one can construct an accurate, cheaply
computable \methodname{}
model for general-output error
$\error = \outbias = \outputfun(\state) - \outputfun(\rightbasiscoarse\redstatecoarse)$ by employing
indicators $\indicators =
\redadjointcoarsetofine^t\modresidual\left(
\rightbasiscoarse\redstatecoarse\right)$ and transformation function $\functionErrorNo=\identityFunction{\RR{}}$ the identity function over $\RR{}$.
\begin{remark}[Uncertainty control for dual-weighted-residual error
indicators]\label{rem:uncertaintyControlDual}
The accuracy of the reduced-order dual solution can be controlled by changing
$\nrbdual$---the dimension of the dual basis $\rightbasisdualcoarse$. In
general, one would expect an increase in $\nrbdual$ to lead to a lower-variance
\methodname{} surrogate at the expense of a higher dimensional dual problem
\eqref{eq:adjointCoarse}. The experiments in Section
\ref{sec:experiments:multipleOutputs} highlight this uncertainty-control
attribute.
\end{remark}

\subsection{Probabilistically rigorous error
bounds}\label{sec:probRigError}
  Clearly, the \methodname{} surrogate does not strictly bound
  the error, even when error bounds are used as indicators. That is, the mean
  probability of overestimation is generally less than one
  \begin{equation}
    \crigor\defeq\frac{1}{|\Para|} \int_{\bmu \in \Para}
    \prob{\quantitysurr(\indicators(\bmu))> \functionError{\error{(\bmu)}}} d\bmu
    < 1.
    \label{eq:probOverestimationPositive}
  \end{equation}

  This frequency of overestimation depends on the probability
  distribution of the random variable $\machine(\indicators)$. Using the
  machine learning methods proposed in the next section, we infer normally
  distributed random variables
  \begin{equation}
    \machine(\indicators) \sim \cN(\nu(\indicators), \overline{\sigma}^2(\indicators))
    \label{eq:gaussianRandomVInferred}
  \end{equation}
  with mean $\nu(\indicators)$ and variance
  $\overline{\sigma}^2(\indicators)$.  If the model is perfectly validated,
	then the mean probability of overestimation is
  $\crigor = 0.5$.  However, knowledge about the distribution of
  the random variable can be used to control the overestimation frequency.
  In particular, the modified surrogate
  \begin{equation}
    {\machine{}}^\crigor(\indicators) \defeq  {\machine{}}(\indicators)
    + \rigorbound(\indicators,\crigor)
    \label{eq:stochstateerrorRigor}
  \end{equation}
  enables \emph{probabilistic  rigor}: it bounds the error with mean
  specified probability $\crigor$ assuming the model is perfectly validated. Here, $\rigorbound$ fulfills
  \begin{equation}
    \prob{X > \rigorbound(\indicators,\crigor)} = \crigor, \qquad \text{ for~}X \sim \cN(0, \overline{\sigma}^2(\indicators)).
  \end{equation}
  This value can be computed as
  \begin{equation}
    \rigorbound(\indicators) = \sqrt{2} \overline{\sigma}(\indicators) \erf^{-1}\rb{2c - 1}
    \label{eq:deltarigor}
  \end{equation}
  where $\erf^{-1}$ is the inverse of the {\em error function}. 

\section{Gaussian processes}
\label{sec:gpbasics}

This section describes the two methods we employ to construct the stochastic
mapping $\quantitysurr:\indicators(\params)\mapsto\errorsurr$:
\begin{\itemname}
  \item {\em Gaussian process kernel regression} (i.e., kriging) 
    \cite{RW06} and
  \item the {\em relevance vector machine} (RVM) \cite{Tipping2001}.
\end{\itemname}
Both methods are examples of supervised learning methods that generate a
stochastic process from a
set of $\romesTraining$ training points for independent variables
$\bx\defeq(\bm{x}_n)_{n=1}^\romesTraining$ and a dependent variable
$\by\defeq(y_n)_{n=1}^\romesTraining$.
Using these training data, the methods generate 
predictions $\machiney(\bm x_m^*;
\theta^{\text{ML}})$, $m=1,\ldots,M$
associated with a set of $M$ prediction
points $\bx^*\defeq (\bm x^*_m)_{m=1}^{M}$. Here,
$\theta^{\text{ML}}$ denotes 
hyperparameters that are inferred using a Bayesian approach; the predictions
are random variables with a multivariate normal distribution.

In the context of \methodname{}, 
the independent variables correspond to error indicators $\bm{x}_n =
\indicators(\params_n)$ with $\params_n\in\paramDomain$,
$n=1,\ldots,\romesTraining$  and the
dependent variable corresponds to the (transformed) reduced-order-model error such that 
$y_n =\functionError{\error(\params_n)}$, $n=1,\ldots,\romesTraining$.
To make this paper as self-contained as possible, the following sections
briefly present and compare the two approaches.

  \subsection{GP kernel method}

  A Gaussian
  process 
  is defined as a collection of random variables such that any finite
  number of them has a joint Gaussian distribution.
  The GP kernel method constructs this Gaussian process via Bayesian inference
  using the training data and a specified kernel function.
  To begin, the prior distribution is set to
\begin{equation} 
\machiney_\mathrm{prior}(\underline \bx)
\sim\mathcal \romesTraining\left(0,\K\left(\underline \bx,\underline \bx\right) + \sigma^2
\bI_{\romesTraining+M}\right)
\end{equation} 
with $\underline\bx \defeq\left(\underline{\bm
x}_i\right)_{i=1}^{\romesTraining+M} =
(\bx,\bx^*)$.
  Here, the GP kernel assumes that the covariance between any two points
  can be described analytically by a kernel $k$ 
with additive noise $\bepsilon \sim \cN(0, \sigma^2\bI_{M+\romesTraining})$  
  such that
  \begin{equation}
    \K(\underline \bx, \underline \bx) = \rb{k(\underline{\bm x}_i, \underline{\bm
    x}_j)}_{1\leq i,j\leq \romesTraining+M}.
    \label{eq:definition:kernel}
  \end{equation}
  In this work, we employ the most commonly used 
  squared-exponential-covariance kernel
  \begin{align}\label{eq:squExpCov}
    k(\bm x_i, \bm x_j) &= \exp\rb{-\frac{\|{\bm x_i-\bm x_j}\|_2^2}{2l^2}},
  \end{align}
  which induces high correlation between geometrically nearby  
  points.  Here, $l\in\RR{}$ is the `width' hyperparameter.

  Assuming the predictions are generated as independent samples from the
	stochastic process,\footnote{Typically in the GP literature, predictions at
	all points $\bx^*$ are generated
	simultaneously as a single sample from the Gaussian process. In this work,
	we treat all predictions as arising from independent samples of the GP.} the GP kernel method then generates predictions for
  each point $\bm x^* \in \bx^*$. These predictions correspond to random
  variables with posterior distributions $\machiney(\bm x^*; \theta) \sim
  \cN(\nu(\bm x^*), \overline{\sigma}^2(\bm x^*))$ with
  \begin{align}
    \nu(\bm x^*) &= \K(\bm x^*, \bx) \rb{\K(\bx, \bx) + \sigma^2
		\bI_\romesTraining}^{-1} \by \\
    \label{eq:twoSigmaContributions}\overline{\sigma}^2(\bm x^*) &= \Sigma(\bm x^*) + \sigma^2\\
    \Sigma(\bm x^*) &= \K(\bm x^*, \bm x^*) - \K(\bm x^*, \bx)\rb{\K(\bx, \bx) + \sigma^2
    \bI_\romesTraining}^{-1}\K(\bx, \bm x^*).
  \end{align}
  More details on the derivation of these expressions can be found in
  Ref.~\cite[ch 2.2]{RW06}. 

  The hyperparameters $\theta\defeq (l^2, \sigma^2)$ can be set to the
  maximum-likelihood values $\theta^{\text{ML}}$ computed as
  the solution to an optimization problem 
  \begin{equation}
    \theta^{\text{ML}} = \arg\max_{\theta} \cL(\theta)
    \label{eq:gp:minimization}
  \end{equation}
with
  the log-likelihood function defined as 
  \begin{equation}
    \cL(l^2,\sigma^2) = -\frac{1}{2} \by^t\rb{\K(\bx, \bx; l^2) + \sigma^2
    \bI_\romesTraining}^{-1}\by - \frac{1}{2} \log\abs{\K(\bx, \bx;
		l^2)+\sigma^2\bI_\romesTraining} -
    \frac{\romesTraining}{2} \log 2\pi.
    \label{eq:definition:gp:loglikelihood}
  \end{equation}
  For details on the derivation of the log likelihood function and 
  problem \eqref{eq:gp:minimization}, we refer to
  Ref.~\cite[ch 5.4]{RW06}.
\begin{remark}\label{rem:uncertaintyContributions}
  The noise component $\sigma^2$ of posterior covariance $\overline{\sigma}^2(\bm x^*)$
  accounts for uncertainty in the assumed GP structure.
  It plays a crucial role for the \methodname{} method: it accounts for the non-uniqueness of the mapping
  $\indicators\mapsto\error$,
  as it is possible for
  $\error(\params_i)\neq\error(\params_j)$ even if
  $\indicators(\params_i) = \indicators(\params_j)$. In particular, this noise
  component represents
  the `information loss' incurred by employing the error indicators in lieu
  of the system inputs as independent variables in the GP.  
  Therefore, this component can be
  interpreted as the \emph{inherent uncertainty in the error due to the
  non-uniqueness of the mapping $\indicators\mapsto\error$}. 
  
  On the other hand, the
  remaining term $\Sigma(\bm x^*)$ of the posterior variance quantifies the uncertainty in the
  mean prediction. This decreases as the number of training points increases.
  Therefore, $\Sigma(\bm x^*)$ can be interpreted as the \emph{uncertainty due to
  a lack of training data.}

  For example, the multifidelity-correction approach employs $\indicators =
  \params$ and therefore should be characterized by $\sigma^2 = 0$, as the
  mapping $\params\mapsto \error$ is unique. However,
  due to the high-dimensional nature of the system-input space $\Para$ in many problems, the
  uncertainty due to lack of training $\Sigma(\bm x^*)$ can be very large unless many
  training points are employed. On the
  other hand, the \methodname{} method aims to significantly reduce
  $\Sigma(\bm x^*)$ by employing a small number of indicators, albeit at the cost of
  a nonzero $\sigma^2$. 
\end{remark}

In light of this remark, we will employ two different types of \methodname{}
models: one that includes the uncertainty due to a lack of training data
(i.e., variance $\overline{\sigma}^2(\bm x^*)$), and one that
neglects this uncertainty (i.e., variance $\sigma^2$).

\subsection{Relevance vector machine (RVM) method}
\label{subsec:rvm}

The RVM \cite{Tipping2001} is based
on a pa\-ra\-meterized discretization of the predictive random variable
\begin{equation}
\machiney(\bm x) = \sum_{\rvmvar=1}^\rvmnum w_\rvmvar \phi_\rvmvar(\bm x) +
\epsilon = \bphi(\bm x)^t
\bw + \epsilon,
\label{eq:definition:rvm:discretization}
\end{equation}
with specified basis functions $\bphi(\bm x)\defeq \left[\phi_1(\bm x)\ \cdots\
\phi_\rvmnum(\bm x)\right]^t\in\RR^{\rvmnum}$, a
corresponding set of random variables $\bw\defeq \left[w_1\ \cdots\
w_\rvmnum\right]^t \in \RR^{\rvmnum}$, with $w_\rvmvar \sim \cN(0, \beta^2_\rvmvar)$ for $\rvmvar=1,\ldots,\rvmnum$ and
noise $\epsilon \sim \cN(0, \sigma^2)$. The hyperparameters ${\bm \beta} =
[\beta_1\ \cdots\ \beta_K]^t
\in \RR^K$ define the prior probability distribution, and are usually chosen by
a likelihood maximization over the training samples. Radial basis functions
\begin{equation}
  \label{eq:rvm:radialBasis}
\phi^{RBF}_\rvmvar(\bm x)=\exp\rb{-\frac1{r^2} \norm{\bar{\bm x}_\rvmvar - \bm
x}_2^2},\quad \rvmvar=1,\ldots,\rvmnum
\end{equation}
constitute the most common choice for basis functions.
For the \methodname{}
method, we often expect a linear relationship  between the indicators and true
errors, likely with a small-magnitude
high-order-polynomial deviation. Therefore, we
also consider Legendre polynomials \cite[Ch.8]{Abramowitz1972}
\begin{equation}
\phi^{Leb}_\rvmvar(\bm x)= P_\rvmvar(\bm x),\quad \rvmvar=1,\ldots,\rvmnum.
\end{equation}
Note that both sets of basis functions are dependent on the training data:
while the centering points $\bar{\bm x}_\rvmvar$, $\rvmvar=1,\ldots, \rvmnum$
in the radial basis functions can be chosen arbitrarily, they are typically
chosen to be equal to the training points. The domain of the Legendre polynomials, on the
other hand, is nominally $[-1, 1]$; therefore the independent variables must
be appropriately scaled to ensure 
the range of training and prediction points $(\bx, \bx^*)$ is included in this
interval.

The RVM method also employs a Bayesian approach to construct the model from
training data. In particular, the vector of hyperparameters ${\bm \beta}$
affects the variance of the Gaussian random variables $\bw$. If these
hyperparameters are computed by a maximum-likelihood or a similar optimization
algorithm, large values for these hyperparameters identify insignificant
components that can be removed.
Therefore, in the \methodname{} context, the RVM can
be used to filter out the least significant error indicators. Apart from this detail, the RVM can be considered a
special case of the GP kernel method with kernel
\begin{equation}
  k(\bm x_i, \bm x_j) = \sum_{\rvmvar=1}^\rvmnum \frac{1}{\beta_\rvmvar}
  \phi_\rvmvar(\bm
  x_i)\phi_\rvmvar(\bm x_j).
\end{equation}
\section{Numerical experiments}
\label{sec:experiments}

  \begin{figure}[t]
    \begin{center}
      \tikzsetnextfilename{thermalblock}
      \begin{tikzpicture}
        \begin{axis}[clip=false, enlargelimits=false, width=4.9cm, xtick=\empty, ytick=\empty]
          \addplot graphics[xmin=0,xmax=1,ymin=0,ymax=1] {figures/thermalblock};
          \node at (axis cs:0.17,0.17) {\large 1};
          \node at (axis cs:0.5,0.17) {\large 2};
          \node at (axis cs:0.82,0.17) {\large 3};
          \node at (axis cs:0.17,0.5) {\large 4};
          \node at (axis cs:0.5,0.5) {\large 5};
          \node at (axis cs:0.82,0.5) {\large 6};
          \node at (axis cs:0.17,0.82) {\large 7};
          \node at (axis cs:0.5,0.82) {\large 8};
          \node at (axis cs:0.82,0.82) {\large 9};
          \draw[green, ultra thick] (axis cs:0,1) -- node[above] {$\Gamma_D$} (axis cs:1,1);
          \draw[red, ultra thick] (axis cs:0,0) -- node[below] {$\Gamma_{N_1}$} (axis cs:1,0);
          \draw[gray, ultra thick] (axis cs:0,0) -- node[left] {$\Gamma_{N_0}$} (axis cs:0,1);
          \draw[gray, ultra thick] (axis cs:1,0) -- (axis cs:1,1);
        \end{axis}
      \end{tikzpicture}
    \end{center}
    \caption{Domain and sample solution $\state(\bmu)$ for the thermal-block
    problem.}\label{fig:thermalblockDomain}
  \end{figure}

This section analyzes the performance of the \methodname{} method on 
Poisson's equation with nine system inputs, using the reduced-basis method
to generate the reduced-order model.
First, \Section
\ref{sec:experiments:thermalblock} introduces the test problem. \Section
\ref{sec:experiments:validation} discusses implementation and validation of
the \methodname{} models. Section
\ref{sec:experiments:comparison:Eldred} compares the \methodname{} method to
the multifidelity-correction approach characterized by employing the model
inputs as error indicators. Section
\ref{sec:experiments:comparison:RBbounds} compares the \methodname{}
stochastic error estimate to the error bound given by the reduced-basis
method.
\Section{\ref{sec:experiments:comparison:surrogates}} compares the performance
of the two machine-learning algorithms: the Gaussian process kernel method and
the relevance
vector machine.
Finally, \Section{\ref{sec:experiments:multipleOutputs}} considers
non-compliant and multiple output functionals, which \methodname{} handles via
dual-weighted-residual error indicators.

  \subsection{Problem setup}
  \label{sec:experiments:thermalblock}

  Consider a finite-element model of heat transport on a square
  domain $\Omega\defeq \cup_{i=1}^9 \Omega_i$ composed of nine parameterized
  materials. The block is cooled along the top boundary to a reference temperature
  of zero, a nonzero heat flux is specified on the bottom boundary, and
  the leftmost boundary is adiabatic. The compliant output functional for
  this problem is defined as the integral over the Neumann domain
  $\Gamma_{N_1}$
  \begin{equation}
    \outputfunCont(\contstate(\bmu)) = \int_{\Gamma_{N_1}} \contstate(\bmu)
    \text{d}x,\quad \bmu \in \Para,
    \label{eq:thermalblock:outputOrig}
  \end{equation}
  where the parameter domain is set to $\Para=[0.1, 10]^9$ and $\contstate$ is
	the continuous representation of the finite-element solution.
  The state variable
  $\contstate(\bmu) \in X= H^1_0\defeq \cb{\contstatedummy \in H^1(\Omega)\;|\;\contstatedummy|_{\Gamma_D} = 0}$
  fulfills the weak form of the parameterized Poisson's equation:
find $\contstate(\bmu) \in X$, such that
  \begin{equation}
    \label{eq:problem:continuous}
    a(\contstate(\bmu), v) = f(v)
    \quad
    \text{for all }v \in X.
  \end{equation}
  Here, the bilinear form $a:X\times X\to X$ and the functional $f:X \to X$
	are defined as
  \begin{equation}
    a(u, v) \defeq \int_{\Omega} \coefficientFunction(x;\bmu) \nabla \contstate(\bmu) \cdot \nabla
		v\, \text{d}x,
    \qquad
    f(v) \defeq \int_{\Gamma_{N_1}} v \,\text{d}x
  \end{equation}
		with boundary conditions
		\begin{gather}
    \nabla \coefficientFunction(x;\bmu) \contstate(\bmu) \cdot n = 0 \quad \text{on }\Gamma_{N_0},
    \qquad \nabla \coefficientFunction(x;\bmu) \contstate(\bmu) \cdot n = 1 \quad \text{on }\Gamma_{N_1}.
    \label{eq:thermalblock:boundaries}
  \end{gather}
  We define the coefficient function $\coefficientFunction:\Omega\times\Para \to \RR$ as
  \begin{equation}
    \coefficientFunction(x;\params) = \sum_{i=1}^9 \mu_{i} \indicator{x}{\Omega_i},
  \end{equation}
  where $\mu_i$ denotes the $i$th component of the parameter vector $\params$,
	and the indicator function $\indicator{x}{A} = 1$ if $x\in A$ and is zero
otherwise.
  \Figure{\ref{fig:thermalblockDomain}} depicts the composition of the domain
  and the location of the boundary conditions.

  By replacing the infinite-dimensional function space $X$ with the
	(finite) $\nrbfine$-dimensional finite-element space $X_{\detI}
  \subset X$ in problem \eqref{eq:problem:continuous},
  one can compute the parameter-dependent state function ${\contfestate}(\bmu) \in
  X_{\detI}$ represented by vectors containing the function's degrees of freedom
	$\state(\params) \in \RR^{\nrbfine}$ (see Section \ref{sec:rb:ellipticProblem}). In
	the experiments, the domain is discretized by triangular finite
	elements, which results in a finite-element
	space $X_{\detI}$ of dimension $\ndof= 10^4$. The high-fidelity output (in
	the notation of Section \ref{sec:problem:full}) is then 
    $\outputfun(\state(\bmu)) \defeq \outputfunCont(\contfestate(\bmu))$, $\bmu \in
		\Para$.

  As described in Section \ref{sec:rb:greedy}, we employ a greedy
	algorithm\footnote{
  All reduced-basis computations are conducted with the
  reduced-basis library \rbmatlab{} (\url{http://www.morepas.org/software/}).
		}
	to generate a reduced-basis space $X_{\redI} \subset X_{\detI}$
  of dimension $\nrb \ll \nrbfine$. The algorithm employs a training set of
	$100$ randomly selected points (i.e., $\card{\Paragreedy} =
	100$
	in Section \ref{sec:rb:greedy}), until the maximum 
  computed error bound in the training set is less than $1$;  it stops after
  $\nrb = 11$ iterations.
	
	Replacing $X_{\detI}$ with $X_{\redI}$ in Eq.~\eqref{eq:problem:continuous} leads to
  reduced state functions ${\contredstate}(\params) \in X_{\redI}$ for all $\bmu
  \in \Para$. As before, these solutions can 
  be represented by vectors $
	\redstate(\params)\defeq\rightbasis\redstatecoords(\params)\in
	\RR^{\nrbfine}$, where
$\rightbasis\in\RR^{\ndof\times \nrb}$ is 
	 the discrete representation of a basis for
	the function space $X_{\redI}$.

	In the following, we analyze two types of error: (i) the energy norm of the
	state-space error
  $
    \enstateerror{} = \vertiii{ \state- \redstate }
    \defeq  a({\contfestate} - {\contredstate},
    {\contfestate} - {\contredstate})
  $
  and (ii) the output error $\outbias = \outputfun({\contfestate}) -
  g({\contredstate})$.  Because the output functional in this case is
	compliant (i.e., $\outputfun = f$ and $a$ is symmetric), the output
	error is always non-negative; see Eq.~\eqref{eq:outboundDefinition} of
	Section \ref{sec:ellipticProblem:errorBounds}.  For more details
	regarding the finite-element discretization, the
	reduced-order-model generation, and error bounds, consult Section
	\ref{sec:rb} of the Supplementary Materials.

  \subsection{\methodname{} implementation and validation\label{sec:experiments:validation}}

  We first compute \methodname{} surrogates for the two errors
	$\enstateerror{}$ and (compliant)
	$\outbias{}$. As proposed in \Section{\ref{sec:errorboundindicators}}, the
	three \methodname{} ingredients we employ are: 1) log-residual-norm error indicators  $\indicators(\params) =
	\log(\residualnorm(\bmu))$, 2) a logarithmic transformation function
	$d=\log$, and 3) both the GP kernel and the RVM supervised
	machine-learning methods.
  To train the surrogates, we compute 
$\enstateerror{(\bmu)}$, $\outbias{(\bmu)}$, and
$\indicators(\bmu)$
   for  $\bmu\in\Pararomes\subset\Para$ with
  $\card{\Pararomes} = 2000$. The first $\romesTraining =100$
  points comprise the \methodname{} training set
	$\{\params_1,\ldots,\params_\romesTraining\}\eqdef\Paratrainromes\subset\Pararomes$  and the following $1900$ points define a validation
  set $\Paravalidate\subset\Pararomes$; note that the validation set was not
	used to construct the error surrogates. Reported results relate to
	statistics computed over this validation set.

  For the kernel method, we employ the squared exponential covariance kernel
	\eqref{eq:squExpCov}. For the RVM method, we choose Legendre polynomials $P_k$ as basis functions, as
	we expect a linear relationship between the indicators and true errors (see
	\Section{\ref{subsec:rvm}}). Because Legendre polynomials are defined on the interval
  $[-1, 1]$, we must transform and scale this domain to span the possible range
  of indicator values. For this purpose, we apply the heuristic of setting
	the
	domain of the polynomials to be 20\% larger than the interval bounded by 
  the smallest and largest indicator values:
  \begin{equation}
  \eb{\indicators_{\min} - 0.1 (\indicators_{\max}-\indicators_{\min}),
    \indicators_{\max} + 0.1 (\indicators_{\max}-\indicators_{\min})},
  \end{equation}
	where
  $\indicators_{\min} = \min_{\bmu\in\Paratrainromes}\indicators(\bmu)$ 
and
$\indicators_{\max} = \max_{\bmu\in\Paratrainromes}\indicators(\bmu)$.
	We include Legendre polynomials of orders $0$ to $4$; however, the RVM
	method typically discards the higher order polynomials due to the near-linear
	relation between indicators and errors.

  \begin{figure}[t]
    \begin{center}
      \tikzsetfigurename{gpviz}
      \begin{tikzpicture}
      {
        \pgfplotsset{
          tick label style={font=\small},
          legend style={font=\small, legend columns=4, legend cell align=center},
          label style={font=\small, text width=4.0cm, align=center},
          xlabel style={font=\small, text width=4.0cm, sloped like x axis, align=center},
          width=5.2cm,}

        \inputaxis{\dirdefault/gaussian_processoutput_index_1_rb_size_18_dual_rb_size_21_v1_tikz.tex}{%
          name=gpviz1, legend to name=leg:gpviz1,%
          title={(i) Kernel method}, ymin=0.1e-2}
        \inputaxis{\dirdefault/gaussian_processoutput_index_1_rb_size_18_dual_rb_size_21_v2_tikz.tex}{%
          at={($(gpviz1.east)+(3cm,0)$)}, anchor=west, legend to name=leg:gpviz2,%
          title={(ii) RVM}, ymin=0.5e-3, ymax=20}

        \inputaxis{res2xnorm/gaussian_processoutput_index_1_rb_size_18_dual_rb_size_21_v1_tikz.tex}{%
          at={($(gpviz1.south)+(0,-0.8cm)$)}, anchor=north,name=gpviz3, legend to name=leg:gpviz3,%
          title={}, ymin=0.1e-2}
        \inputaxis{res2xnorm/gaussian_processoutput_index_1_rb_size_18_dual_rb_size_21_v2_tikz.tex}{%
          at={($(gpviz3.east)+(3cm,0)$)}, anchor=west, legend to name=leg:gpviz4,%
          title={}, ymin=0.5e-3, ymax=20}
        \node at ($(gpviz3.south east)+(1.5cm, -1cm)$) [anchor=north] {\pgfplotslegendfromname{leg:gpviz3}};
      }
      \end{tikzpicture}
      \caption{Visualization of \methodname{} surrogates ($\error = \enstateerror{}$ and $\stateerror{}_X$, $\indicators = \log{\residualnorm}$, $\functionErrorNo =
\log$), computed
      using $\romesTraining =100$ training points and the (i) GP kernel method and (ii) RVM.
			}\label{fig:gpvisualization}
    \end{center}
  \end{figure}

  \Figure{\ref{fig:gpvisualization}} depicts the \methodname{} surrogate
  $\stochastic{\enstateerror{}}$ generated by both machine-learning methods
  using all $100$ training points. 
	For comparison, we also create \methodname{} surrogates
	$\stochastic{\stateerror{}_X}$ for errors in the parameter-independent norm
	$\norm{\cdot}_X$ of the state space $X=H^1_0$. 
	In addition to the
  expected mean of the inferred surrogate, the figure displays two
  $95\%$-confidence intervals for the prediction (see Remark \ref{rem:uncertaintyContributions}):
  \begin{\itemname}
    \item The darker shaded interval  corresponds to the confidence interval
      arising from the inherent uncertainty in the error due to the
      non-uniqueness of the mapping $\indicators\mapsto\enstateerror{}$, i.e.,
      the inferred variance $\sigma^2$ of Eq.~\eqref{eq:twoSigmaContributions}.
    \item The lighter shaded interval also includes the `uncertainty in the
      mean' due to a lack of training data, i.e., $\Sigma$ of
      Eq.~\eqref{eq:twoSigmaContributions}. With an increasing number of
      training points, this area should be indistinguishable from the darker
      one.
  \end{\itemname}

	All \methodname{} models find a linear trend between the indicators and the
	errors, where the variance is slightly larger for the parameter-independent
	norm. This larger variance can be attributed to the larger range of the
	coercivity constants the parameter-independent norm (see Section
	\ref{sec:ellipticProblem:errorBounds}). For this example, however, both
	\methodname{} are functional. In the following examples, we focus on the
	energy norm only.

  Note that the `uncertainty in the mean' is dominant for the RVM
  surrogate. This can be explained as follows: the high-order polynomials have
  values close to zero near the mean of the data. As such, the
  training data are not very informative for the coefficients of
  these polynomials. This results in a large inferred variance for those
  coefficients.
  \Section{\ref{sec:experiments:comparison:surrogates}} further
  compares the two machine-learning methods; due to its superior performance,
	we now proceed with the kernel method.

We now assess the validity of the Gaussian-process assumptions underlying the
\methodname{} surrogates $\stochastic{\enstateerror{}}$ and
$\stochastic{\outbias}$, i.e., Condition \ref{cond:validated} of
\Section{\ref{sec:statModel}}.
From the discussion in Remark~\ref{rem:uncertaintyContributions}, we know
if the underlying GP model form is correct, then as the number of training
  points increases, the uncertainty about the mean decreases and the set
  $\cb{D(\bmu)\ |\ \bmu \in \Paravalidate}$ with 
  \begin{equation}
    D(\bmu) \defeq  \functionErrorNo\rb{\enstateerror{(\bmu)}} - \expected{\functionErrorNo\rb{\stochastic{\enstateerror{}}(\indicators(\bmu))}}
    = \functionErrorNo\rb{\enstateerror{(\bmu)}} - \nu\rb{\indicators(\bmu)}
    \label{eq:binplotdeviation}
  \end{equation}
should behave like samples from  the distribution $\cN(0, \sigma^2)$.
  \Figure{\ref{fig:histogramsCompare}} reports this validation test and
	verifies that this condition does indeed hold for a sufficiently large number of training points.

  \begin{figure}[t]
    \centering
      \tikzsetfigurename{histograms_compare}
      \begin{tikzpicture}
        {
        \pgfplotsset{
          tick label style={font=\footnotesize},
          label style={font=\footnotesize},
          legend style={font=\small, legend columns=-1, legend cell align=center},
          width=5cm,
          ylabel style={text width=3cm, align=center}}
          \inputaxis{\dirdefault/histogram_10_bs_\rbmedium_v\vardefault_tikz.tex}{%
            name=hist10, legend to name=leg:hist10, title={$\romesTraining =10$}}
          \inputaxis{\dirdefault/histogram_35_bs_\rbmedium_v\vardefault_tikz.tex}{%
            name=hist25, at={($(hist10.east)+(2cm,0)$)}, anchor=west,%
            legend to name=leg:hist25,%
            title={$\romesTraining =35$}}
          \inputaxis{\dirdefault/histogram_65_bs_\rbmedium_v\vardefault_tikz.tex}{%
            name=hist50, at={($(hist25.east)+(2cm,0cm)$)}, anchor=west,%
            legend to name=leg:hist50,%
            title={$\romesTraining =65$}}
          \inputaxis{\dirdefault/histogram_95_bs_\rbmedium_v\vardefault_tikz.tex}{%
            name=hist95, at={($(hist10.south)-(0cm,2cm)$)}, anchor=north,%
            legend to name=leg:hist75,%
            title={$\romesTraining =95$}}
          \node at ($(hist95.east)+(2cm,0cm)$) [anchor=west] {%
             \small%
             \input{\imagefile{res2energy/confidence_intervals_bs_\rbmedium_v\vardefault_table_tikz.tex}} };
          \node at ($(hist95.south east)+(4.0cm, -1cm)$) [anchor=north] {\pgfplotslegendfromname{leg:hist75}};
        }
      \end{tikzpicture}
    \caption{
			Gaussian-process validation for the \methodname{} surrogate (GP kernel,
			$\error = \enstateerror{}$, $\indicators = \log{\residualnorm}$, $\functionErrorNo =
\log$) with a varying number of training points $\romesTraining $. 
   The histogram corresponds to samples of $D(\bmu)$
	 and the red curve depicts the probability
	 distribution function $\mathcal N(0,\sigma^2)$.  
    The table reports how often the actual error lies in the inferred confidence
    intervals.}\label{fig:histogramsCompare}
  \end{figure}
  Further, we can validate the inferred confidence intervals
  as proposed in Eq.~\eqref{eq:ciValidate}. The table within
	\Figure{\ref{fig:histogramsCompare}} reports $\cobs{\ciVar}$ (see Eq.~\eqref{eq:validatedModelA}), which represents the frequency of observed predictions in the validation
set that lie within the
inferred confidence interval $\confidenceIntervalOne{\ciVar}$.
We declare the \methodname{} model to be validated, as $\cobs{\ciVar}\approx
\ciVar$
for several values of $\ciVar$ as the number of training points increases.

	The results for the \methodname{} surrogate $\stochastic{\outbias{}}$ are
	very similar to those presented in \Figure{\ref{fig:histogramsCompare}} and
	will be further discussed in
	\Section{\ref{sec:experiments:comparison:Eldred}}.
	Note that the inferred Gaussian process is well-converged with a
	moderately sized training set consisting of only $\romesTraining=35$ points.

\subsection{Output error: comparison with multifidelity correction}
\label{sec:experiments:comparison:Eldred}
As discussed in \Section{\ref{sec:romes:stochasticMapping}},
multifidelity-correction methods construct a surrogate $\stochastic{\outbiasMF}$ of the output error
using the system inputs as error-surrogate inputs, i.e., 
$\delta = \outbias$, $\indicators = \params$, and $\functionErrorNo =
\identityFunction{\RR{}}$. In this section, we construct this multifidelity
correction surrogate using the same GP kernel method as \methodname{}.
Ref.~\cite{NE12} demonstrated that this error surrogate
fails to improve the `corrected output' when the low-fidelity model
corresponds to a reduced-order model. We now verify this
result and show that---in contrast to the multifidelity correction
approach---the \methodname{} surrogate $\stochastic{\outbias}$ constructed
via the GP kernel method with 
$\delta = \outbias$, $\indicators = \log{\residualnorm}$, and $\functionErrorNo =
\log$
yields impressive results:
on average, the output `corrected' by the \methodname{} surrogate reduces the
error by an order of magnitude, and the Gaussian-process assumptions can be
validated. The validation quality improves as the number of training points 
increases, but a moderately sized set of only $\romesTraining=20$ training points leads to a
converged surrogate.

\begin{figure}[thp]
  \begin{center}
    \tikzsetfigurename{correlation_output_bias}
    \begin{tikzpicture}
    {
      \pgfplotsset{
        tick label style={font=\small},
        label style={font=\small, text width=4.0cm, align=center},
        legend style={font=\small},
        width=6cm,
      }

      \inputaxis{res2bias/correlation_1d_bs_\rbmedium_v1_tikz.tex}{%
        legend to name=leg:biascorrelation1, name=biascorrelation1,%
        legend entries={$(r;\outbias{})$,$(\outbound;\outbias{})$},%
        title={(i) Output error v.\ \methodname{} indicators}}

      \inputaxis{param2logbias/correlation_1d_bs_\rbmedium_v1_tikz.tex}{%
        at={($(biascorrelation1.east)+(3cm,0)$)}, anchor=west,%
        legend to name=leg:biascorrelation2, name=biascorrelation2,%
        legend entries={$(\mu_1;\outbias{(\bmu)})$,$(\mu_2;\outbias{(\bmu)})$},%
        title={(ii) Output error v.\ system inputs}}

      \node at ($(biascorrelation1.south east)-(0.1cm,-0.1cm)$) [anchor=south east] {\pgfplotslegendfromname{leg:biascorrelation1}};

      \node at ($(biascorrelation2.south east)-(0.1cm,-0.1cm)$) [anchor=south east] {\pgfplotslegendfromname{leg:biascorrelation2}};
     }
    \end{tikzpicture}
  \end{center}
  \caption{Relationship between (i) \methodname{} error indicators
	and the compliant-output error
  and (ii) the first two parameter components and the (compliant) output
	error, visualized by evaluation of 200 random sample points in the input
	space.
	Clearly,
  the observed structure in the former relationship is more amenable to
	constructing a Gaussian process.}\label{fig:biascorrelations}
\end{figure}
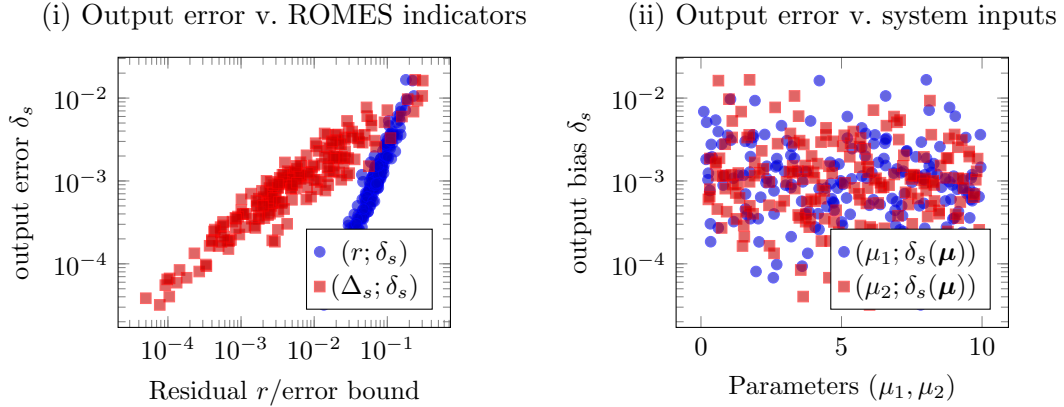

The reason multifidelity correction fails for most reduced-order models is
twofold. First, the 
mapping
$\bmu \mapsto \outbias$ can be highly
oscillatory in the input space. This behavior arises from the fact the the
reduced-order model error is zero at the (greedily-chosen) ROM training points but
grows (and can grow quickly) away from these points. Such complex behavior
requires a large number of error-surrogate training points to accurately
capture. In addition, the number of system inputs is often large (in this case
nine); this introduces curse-of-dimensionality difficulties in modeling the error. \Figure{\ref{fig:biascorrelations}(ii)} visualizes this problem. The depicted
mapping between the first two parameter components $\mu_1, \mu_2$ and the
output error  $\outbias(\bmu)$ displays no structured behavior.
As a result, there is no measurable improvement of the
corrected output $\redout{} + \stochastic{\outbiasMF{}}$ over the evaluation
of the ROM output $\redout$ alone.

In order to quantify the performance of the error surrogates, we introduce a normalized
{\em expected improvement}
\begin{equation}
  I(\stochastic{\delta}, \bmu) \defeq \left|\frac{\outbias{(\bmu)} -
  \mode{\stochastic{\delta{}}(\indicators(\bmu))}}{\outbias{(\bmu)}}\right|.
  \label{eq:improvementValueMF}
\end{equation}
If this value is less than one, then the exected corrected output $\redout{} +
\stochastic{\error{}}$ is more accurate than the ROM output
$\redout{}$ itself for point $\bmu\in\Para$, i.e., the additive error surrogate \emph{improves} the prediction
of the ROM. On the other hand, values above one indicate that the error
surrogate \emph{worsens} the ROM prediction.

\Figure{\ref{fig:improvementComparisonEldred}} reports the mean, median,
standard deviations, and extrema for the expected improvement 
\eqref{eq:improvementValueMF} evaluated for all validation points 
$\Paravalidate$ and a varying number of training points. 
Here, we also compare with the performance of
the error surrogate $\errorSurrUni$, which is defined as a uniform
distribution on the interval
$\left[\outlowerbound{(\bmu)},\outbound{(\bmu)}\right]$, where
$\outlowerbound{(\bmu)}$  and $\outbound{(\bmu)}$ are the 
the lower and upper bounds for the output error, respectively. Note that
$\errorSurrUni$ does not require training data, as it is based purely on
error bounds.

The expected improvement for the \methodname{} output-error surrogate
$I(\stochastic{\outbias{}}, \bmu)$ as depicted in
\Figure{\ref{fig:improvementComparisonEldred}(i)} is approximately $0.2$ on
average, which constitutes an improvement of nearly an order of magnitude.
Further, the maximum expected improvement almost always remains below $1$; this
implies that the corrected ROM output is almost always more accurate than
the ROM output alone.

On the other hand, the expected improvement generated by
the error surrogate $\errorSurrUni$ is always greater than one, which means
that its correction always increases the error. This arises from the fact that
the center of the interval
$\left[\outlowerbound{(\bmu)},\outbound{(\bmu)}\right]$ is a poor
approximation for the true error.

In addition, \Figure{\ref{fig:improvementComparisonEldred}(ii)} shows that the
expected improvement produced by the multifidelity-correction surrogate
$I\rb{\stochastic{\outbiasMF}, \bmu}$ is often far greater than one.  This
shows that the multifidelity-correction approach is not well suited for this
problem. Presumably, with (far) more training points, these results would
improve.

  \begin{figure}[tp]
    \begin{center}
      \tikzsetfigurename{improvementComparisonEldred}
      \begin{tikzpicture}
      {
        \pgfplotsset{
          tick label style={font=\small},
          legend style={font=\small, legend columns=4, legend cell align=center},
          label style={font=\small, text width=4.7cm, align=center},
          title style={text width=4.6cm, align=center},
          width=6cm}
        \inputaxis{res2bias/improvement_mode_bs_\rbmedium_v\vardefault_tikz.tex}{%
          legend to name=leg:res2biasImprovement, name=leg2biasImprovement,%
          ymin=0.001, title={(i) \methodname{}}}
        \inputaxis{param2logbias/improvement_bs_\rbmedium_v1_tikz.tex}{%
          at={($(leg2biasImprovement.east)+(3cm,0)$)}, anchor=west,%
          restrict y to domain=-100:100,%
          ymin=0.01,
          legend to name=leg:param2biasImprovement,%
          title={(ii) Multifidelity correction}}
        \node at ($(leg2biasImprovement.south east)+(1.5cm, -1.1cm)$) [anchor=north] {\pgfplotslegendfromname{leg:res2biasImprovement}};
      }
      \end{tikzpicture}

    \end{center}
    \caption{Expected improvement $I(\stochastic{\delta}, \bmu)$ for a varying
		number of training points $\romesTraining $: (i) \methodname{} (GP kernel, compliant $\delta =
		\outbias$, $\indicators = \log{\residualnorm}$, $\functionErrorNo =
\log$) with uniform distribution based on reduced-basis error bounds
$\errorSurrUni$ and (ii) multifidelity correction (GP kernel, compliant $\delta = \outbias$, $\indicators = \params$, and $\functionErrorNo =
\identityFunction{\RR{}}$). ($1$: no
    improvement, $>1$: error worsened, ${}<1$: error improved).}\label{fig:improvementComparisonEldred}
  \end{figure}

  \begin{figure}[thp]
    \centering
      \tikzsetfigurename{histogramsCompareEldred}
      \begin{tikzpicture}
        {
        \input{\imagefile{confidence_intervals_res2bias_param2logbias_bs_\rbmedium_v\vardefault.tex}}
        \pgfplotsset{
          tick label style={font=\footnotesize},
          label style={font=\footnotesize, text width=3.6cm, align=center},
          legend style={font=\small, legend columns=-1, legend cell align=center},
          title style={text width=4.5cm, align=center},
          width=5cm}
          \inputaxis{res2bias/histogram_100_bs_\rbmedium_v\vardefault_tikz.tex}{%
            name=histres2bias, legend to name=leg:histres2bias, title={(i)
						\methodname{}}}
          \inputaxis{param2logbias/histogram_100_bs_\rbmedium_v1_tikz.tex}{%
            name=histparam2bias, at={($(histres2bias.east)+(1.5cm,0)$)}, anchor=west,%
            legend to name=leg:histparam2bias,%
            title={(ii) Multifidelity correction}}
          \node at ($(histparam2bias.east)+(0.5cm,0cm)$) [anchor=west] {%
             \small%
             \cintsrestobiasparamtobias{100}{(i)}{(ii)} };
          \node at ($(histparam2bias.south)+(0.0cm, -1cm)$) [anchor=north] {\pgfplotslegendfromname{leg:histparam2bias}};
        }
      \end{tikzpicture}

    \caption{
			Gaussian-process validation for the \methodname{} surrogate (GP kernel,
			compliant $\delta = \outbias$, $\indicators = \log{\residualnorm}$, $\functionErrorNo =
\log$) and
			multifidelity-correction surrogate (GP kernel, compliant $\delta = \outbias$, $\indicators = \params$, and $\functionErrorNo =
\identityFunction{\RR{}}$) using $\romesTraining =100$ training points
   The histogram corresponds to samples of $D(\bmu)$
	  and the red curve depicts the probability
	 distribution function $\mathcal N(0,\sigma^2)$.  
  The table reports how often the
    actual error lies in the inferred confidence intervals. Clearly, this
		validation test fails for the multilfidelity-correction surrogate.
    \label{fig:histogramsCompareEldred}
    }
  \end{figure}

Again, we can validate the Gaussian-process assumptions underlying the error
surrogates. For $\romesTraining =100$ training points, \Figure{\ref{fig:histogramsCompareEldred}} compares a histogram
of deviation of the true error from the surrogate mean
to the inferred probability density function.
The associated table reports how often the
validation data lie in inferred confidence intervals. We observe that
the confidence intervals are valid for the \methodname{} surrogate, but are
not for the multifidelity-correction surrogate, as the bins
do not align with the inferred distribution. \Figure{\ref{fig:cintsCompareEldred}} depicts the convergence of these
confidence-interval validation metrics as the number of training points
increases. As expected (see Remark~\ref{rem:uncertaintyContributions})
the \methodname{} observed confidence intervals more closely align with the confidence
intervals arising from the \emph{inherent uncertainty} (i.e., 
$\sigma^2$) as the number of training points increases, as this effectively
decreases the \emph{uncertainty due to a lack of training}.
In addition, only a moderate number of training points 
(around $20$) is required to generate a reasonably converged \methodname{}
surrogate. On the other hand the multifidelity-correction surrogate exhibits
no such convergence when fewer than 100 training points are used.

  \begin{figure}[htp]
    \centering
      \tikzsetfigurename{cintsCompareEldred}
      \begin{tikzpicture}
        {
        \pgfplotsset{
          tick label style={font=\footnotesize},
          label style={font=\footnotesize, text width=4.5cm, align=center},
          legend style={font=\small, legend columns=-1, legend cell align=center},
          title style={text width=4.6cm, align=center},
          width=6cm}
          \inputaxis{res2bias/confidence_intervals_plot_bs_\rbmedium_v\vardefault_tikz.tex}{%
            name=cintallres2bias, legend to name=leg:cintallres2bias,
						title={(i) \methodname{}}}
          \inputaxis{param2logbias/confidence_intervals_plot_bs_\rbmedium_v1_tikz.tex}{%
            name=cintallparam2bias, at={($(cintallres2bias.east)+(2cm,0)$)}, anchor=west,%
            legend to name=leg:cintallparam2bias,%
            title={(ii) Multifidelity correction}}
          \node at ($(cintallres2bias.south east)+(1cm, -1cm)$) [anchor=north] {\pgfplotslegendfromname{leg:cintallres2bias}};
        }
      \end{tikzpicture}

    \caption{Gaussian-process validation for the \methodname{} surrogate (GP
		kernel, compliant $\delta = \outbias$, $\indicators = \log{\residualnorm}$, $\functionErrorNo =
\log$) and
			multifidelity-correction surrogate (GP kernel, compliant $\delta = \outbias$, $\indicators = \params$, and $\functionErrorNo =
\identityFunction{\RR{}}$) and a varying number of training points $\romesTraining $.
		The plots depict how often the
    actual error lies in the inferred confidence intervals.
    }\label{fig:cintsCompareEldred}
  \end{figure}

  \subsection{Reduced-basis error bounds}
  \label{sec:experiments:comparison:RBbounds}
  In this section, we compare the reduced-basis error bound
  $\statebound^{\bmu}$ \eqref{eq:stateboundDefinitionEnergy} with the
	\emph{probabilistically rigorous} \methodname{} surrogates
  $\stochastic{\enstateerror{}}^c$  \eqref{eq:effectivityBoundDefinition} with 
	rigor values of $\crigor=0.5$ and $\crigor = 0.9$ as introduced Section
	\ref{sec:probRigError}.\footnote{Note that $\crigor=0.5$ implies no modification to the
	original \methodname{}
	surrogate, as $ \stochastic{\enstateerror{}}^{0.5}=\stochastic{\enstateerror{}}$ (see
	Eqs.~\eqref{eq:stochstateerrorRigor}--\eqref{eq:deltarigor}).}
  The \methodname{} surrogate is
	constructed with the GP kernel method and ingredients $\delta =
	\enstateerror{}$, $\indicators = \log{\residualnorm}$, and $\functionErrorNo
  = \log$. 
	As
	discussed
  in \Section{\ref{sec:problem:errors}} the error-bound effectivity 
  \eqref{eq:effectivityBoundDefinition} is important to quantify the
  performance of these bounds; a value of 1 is optimal, as it implies no
	over-estimation.
As the probabilistically rigorous \methodname{} surrogates
$\stochastic{\enstateerror{}}^\crigor$ are stochastic processes, we can measure
	their (most common) effectivity as
			\begin{equation}
				\eta(\crigor,\bmu) \defeq
				{\frac{\mode{\stochastic{\enstateerror{}}^\crigor\rb{\indicators(\bmu)}}}{\enstateerror{(\bmu)}}}
				\label{eq:stochasticEta}.
			\end{equation}

		\begin{figure}[t]
			\begin{center}
				\tikzsetfigurename{effectivityRBbound}
				\begin{tikzpicture}
				{
					\pgfplotsset{
						tick label style={font=\small},
						legend style={font=\small, legend columns=1, legend cell align=center},
						label style={font=\small, text width=4.7cm, align=center},
						title style={text width=4.5cm, align=center},
						width=5.2cm}
					\inputaxis{\dirdefault/effectivity_mode_50_with_delta_bs_\rbmedium_v1_tikz.tex}{%
						legend to name=leg:effectivityRBbound, name=effectivityRBbound, ymode=linear,%
						title={$50\%$-rigorous estimate}}
						\quad
					\inputaxis{\dirdefault/effectivity_90_with_delta_bs_\rbmedium_v1_tikz.tex}{%
						at={($(effectivityRBbound.east)+(2cm,0)$)}, anchor=west,%
						title={$90\%$-rigorous estimate},%
						legend to name=leg:effectivityRBbound90, name=effectivityRBbound90, ymode=linear}
					\inputaxis{\dirdefault/rigor_50_bs_\rbmedium_v1_tikz.tex}{%
						at={($(effectivityRBbound.south)+(0, -2cm)$)}, anchor=north,%
						name=rigor50, ymode=linear}
						\quad
					\inputaxis{\dirdefault/rigor_90_bs_\rbmedium_v1_tikz.tex}{%
						at={($(effectivityRBbound90.south)+(0, -2cm)$)}, anchor=north,%
						name=rigor90, ymode=linear}
					\node at ($(effectivityRBbound90.east)+(0.5cm, 0cm)$) [anchor=west] {\pgfplotslegendfromname{leg:effectivityRBbound}};
				}
				\end{tikzpicture}
			\end{center}
			\caption{Validation of the probabilistically rigorous \methodname{}
			surrogates $\stochastic{\enstateerror{}}^\crigor$ 
		(GP kernel, $\delta = \enstateerror{}$, $\indicators = \log{\residualnorm}$, $\functionErrorNo = \log$)	
			and comparison with RB error upper bound $\statebound^{\bmu}$ and uniform distribution based on reduced-basis error bounds
	$\errorSurrUni$.
			The top plots compare statistics of the effectivities $\eta(\crigor,\bmu)$
			with $\crigor=0.5$ and
			$\crigor=0.9$  of the probabilistically rigorous
			\methodname{} surrogates with the RB error-bound surrogates. The bottom
			plots compare the frequency of error overestimation 
			$\cvalidation$ with the desired value $\crigor$ (red
			line).
			}\label{fig:effectivityRBbound}
		\end{figure}

		The top plots of \Figure{\ref{fig:effectivityRBbound}} report the
		mean, median, standard deviation, and extrema of the 
  effectivities $\eta(0.5,\bmu)$ and $\eta(0.9,\bmu)$ for all
  validation points $\bmu\in\Paravalidate$. Again, we compare with
	$\errorSurrUni$, which is a uniform distribution on an interval whose
	endpoints correspond to the
	lower and upper bounds for the error $\enstateerror{(\bmu)}$.
  We also compare with the corresponding statistics for the
	effectivity of the RB error bound $\statebound^{\bmu}$.
  The lower bound for the coercivity constant that is needed in the RB error
  bound $\statebound^{\bmu}$ is chosen as the minimum over all parameter
  components
  $
    \alpha_{\text{LB}}(\bmu)=\min_{i\in\{1,\ldots,9\}} \mu_i.
  $
  This simple choice is effective because the example is affinely parameter
	dependent and linear \cite[Ch. 4.2]{PR07}.

We observe that
  the \methodname{} surrogate yields better results than both the error bound
	$\statebound^{\bmu}$ (which produces effectivities
	roughly between one and eight)
	and the uniform
	distribution $\errorSurrUni$ (which produces mode effectivities 
	roughly between one and four).
  The $50\%$-rigorous \methodname{} surrogate has an almost
  perfect mean effectivity of $1$ as desired. The $90\%$-rigorous surrogate
	has a higher mean effectivity as expected; however, it is only slightly
	higher.
  Furthermore, the effectivities of the \methodname{} surrogate exhibit a much
	smaller variance\footnote{The higher variance apparent between 45 and 53 training points can be
	explained by the fact that the minimization
  algorithm for the log--likelihood function stops after it exceeds the maximum
	number of iterations.} than both $\statebound^{\bmu}$ and
	$\errorSurrUni$.  Finally, a moderate number (around 20) of training
	samples is sufficient to obtain well-converged surrogates.

 The bottom plots of \Figure{\ref{fig:effectivityRBbound}} report
the frequency of error overestimation 

\begin{equation}
  \cvalidation\defeq\frac{
    \card{\{\bmu\in\Paravalidate\
    |\
    \median{\quantitysurr\rb{\indicators(\bmu)}} >
    \functionError{\error{(\bmu)}}\}}
  }{
    \card{\Paravalidate}
  }
\end{equation}
for the probabilistically rigorous \methodname{} surrogates (i.e., $\quantitysurr =
\stochastic{\enstateerror{}}^\crigor$) as the number of
training points increases to show that $\cvalidation\approx\crigor$, which
validates the rate of error overestimation (see
Eq.~\eqref{eq:probOverestimationPositive}).  Note that the overestimation
frequency $\cvalidation$ converges to its predicted value $\crigor$, which
demonstrates that the rigor of the \methodname{} estimates can in fact be
controlled.

  \subsection{Comparison of machine-learning algorithms}
  \label{sec:experiments:comparison:surrogates}
	This section compares in detail the \methodname{} surrogates generated using the two machine-learning methods introduced in
	\Section{\ref{sec:gpbasics}}. Recall that
	\Figure{\ref{fig:gpvisualization}} visualizes both surrogates. We observe
	that both methods work well overall and
	generate well-converged surrogates with a modest number of training samples.
	As previously mentioned, the GP kernel leads to a smaller inferred 
	variance due to more accurate and localized estimates of the mean.
	The RVM is characterized by global basis functions that
	preclude it from accurately resolving localized features of the mean and
	lead to large uncertainty about the mean (see \Figure{\ref{fig:gpvisualization}}). On the other
	hand, the confidence intervals computed with the RVM
	are (slightly) better validated.

  \begin{figure}[thp]
    \centering
      \tikzsetfigurename{histogramsCompareRVM}
      \begin{tikzpicture}
        {
        \input{\imagefile{histogram_res2energy_kernel_gp_bs_\rbmedium.tex}}
        \pgfplotsset{
          tick label style={font=\footnotesize},
          label style={font=\footnotesize, text width=3.6cm, align=center},
          legend style={font=\small, legend columns=-1, legend cell align=center},
          title style={text width=4.5cm, align=center},
          width=4.5cm}
          \inputaxis{\dirdefault/histogram_80_bs_18_v1_tikz.tex}{%
            name=histreskernel, legend to name=leg:histreskernel, title={(i) Kernel method}}
          \inputaxis{\dirdefault/histogram_80_bs_\rbmedium_v2_tikz.tex}{%
            name=histresrvm, at={($(histreskernel.east)+(1.5cm,0)$)}, anchor=west,%
            legend to name=leg:histresrvm,%
            title={(ii) RVM}}
          \node at ($(histresrvm.north east)+(0.5cm,2em)$) [anchor=north west] {%
             \small%
             \cintsrestoenergykernelgp{80}{(i)}{(ii)}};
          \node at ($(histresrvm.south west)-(0.5cm, 1.4cm)$) [anchor=north] {\pgfplotslegendfromname{leg:histresrvm}};
        }
      \end{tikzpicture}
    \caption{
			Gaussian-process validation for the \methodname{} surrogate using both
			the GP kernel method and the RVM ($\error = \enstateerror{}$, $\indicators = \log{\residualnorm}$, $\functionErrorNo =
\log$) using $\romesTraining  = 80$ training points.
   The histogram corresponds to samples of $D(\bmu)$
	  and the red curve depicts the probability
	 distribution function $\mathcal N(0,\sigma^2)$.
    The table reports how often the actual error lies in the inferred confidence
    intervals.}\label{fig:histogramsCompareRVM}
  \end{figure}

	\Figure{\ref{fig:histogramsCompareRVM}(i) and (ii)} report the validation
	test, i.e., the frequency of deviations from the inferred mean $D(\bmu)$
	\eqref{eq:binplotdeviation} with a training sample containing $\romesTraining =80$
	training points. We observe a smaller inferred variance $\sigma^2$ for the
	GP kernel method, which implies that the mean is identified more
	accurately. In both cases, the validation samples align well with the
	probability density function of the inferred distribution $\cN(0,
	\sigma^2)$.

  The confidence intervals of this inferred distribution can be validated, and
  they turn out to be (slightly) more realistic for the RVM method. The table
	within \Figure{\ref{fig:histogramsCompareRVM}} shows that the kernel
  method results are usually optimistic, i.e., the actual confidence intervals
  are smaller than predicted.  This effect can be corrected, however, by adding
  $\Sigma(\bm x^*)$ as an indicator of the uncertainty of the mean
  as discussed in Remark \ref{rem:uncertaintyContributions}. However, doing so
	for the RVM prediction yields extremely wide confidence intervals due to
	the significant uncertainty about the RVM mean (see
	\Figure{\ref{fig:gpvisualization}}).

  \begin{figure}[thp]
    \begin{center}
      \tikzsetfigurename{effectivityComparisonRVM}
      \begin{tikzpicture}
      {
        \pgfplotsset{
          tick label style={font=\small},
          legend style={font=\small, legend columns=4, legend cell align=center},
          label style={font=\small, text width=4.5cm, align=center},
          title style={text width=4.5cm, align=center},
          width=6cm}
        \inputaxis{\dirdefault/effectivity_mode_50_bs_\rbmedium_v1_tikz.tex}{%
          legend to name=leg:kernelEffectivity, name=kernelEffectivity,%
          title={(i) Kernel method}}
        \inputaxis{\dirdefault/effectivity_mode_50_bs_\rbmedium_v2_tikz.tex}{%
          at={($(kernelEffectivity.east)+(3cm,0)$)}, anchor=west,%
          legend to name=leg:rvmEffectivity,%
          title={(ii) RVM}}
        \node at ($(kernelEffectivity.south east)+(1.5cm, -1cm)$) [anchor=north] {\pgfplotslegendfromname{leg:rvmEffectivity}};
      }
      \end{tikzpicture}
    \end{center}
    \caption{Comparison of the effectivity $\eta(0.5,\bmu)$ of \methodname{}
      surrogates ($\error = \enstateerror{}$, $\indicators =
      \log{\residualnorm}$, $\functionErrorNo = \log$) with the GP computed via
      (i) the GP kernel and (ii) the relevance vector machine
      method.  }\label{fig:effectivityComparisonRVM}
  \end{figure}

  As the inferred variance is larger for the relevance vector machine, this
  also affects the performance of effectivity and improvement measures for the
  error surrogates. \Figure{\ref{fig:effectivityComparisonRVM}} depicts
  statistics of 
  $\eta(0.5,\bmu)$ computed with both the methods,
  and we observe that all statistical measures are significantly better for the
  kernel method estimate.

	We conclude that while both methods produce feasible \methodname{}
	surrogates, the GP kernel method produces consistently better results. In
	particular, the lower inferred variance implies that a lower amount of
	\emph{epistemic uncertainty} is introduced by the error surrogate (See
	condition \ref{cond:lowvar} from \Section{\ref{sec:statModel}}). This is
	critically important for many UQ tasks.
	Therefore, we recommend the GP kernel method to construct
	\methodname{} surrogates.

  \subsubsection{Dependence on reduced-basis size}
  \label{sec:experiments:largerRB}
	To assess the generalizability of the \methodname{} method, we apply it to
	a ROM of higher fidelity, i.e., larger $\nrb$.
	We construct two \methodname{} surrogates: one for the state-space error
	$\stochastic{\enstateerror{}}$ (GP kernel, $\error = \enstateerror{}$,
	$\indicators = \log{\residualnorm}$, $\functionErrorNo = \log$) and one
	for the compliant-output error $\stochastic{\outbias}$ (GP kernel, $\delta = \outbias$, $\indicators = \log{\residualnorm}$, $\functionErrorNo =
\log$). 

	We increase the reduced-basis size by decreasing the maximum error over the
	training set from $1.0$ to $1.0\times
	10^{-3}$ and applying the greedy
	method. The resulting reduced-basis dimension is $\nrb = 62$.
	\Figure{\ref{fig:gpvisualizationBig}} reports the error data and
	\methodname{} surrogates. Comparing the leftmost plot of \Figure{\ref{fig:gpvisualizationBig}}
	with \Figure{\ref{fig:gpvisualization}(i)} and the rightmost plot of
	\Figure{\ref{fig:gpvisualizationBig}} with
	\Figure{\ref{fig:biascorrelations}(i)} reveals that while the errors are
	several orders of magnitude smaller for the current experiments, the data
	(and the resulting \methodname{} surrogates)
	exhibit roughly the same structure as before.

\begin{figure}[htp]
  \begin{center}
    \tikzsetfigurename{gpvizBig}
    \begin{tikzpicture}
    {
      \pgfplotsset{
        tick label style={font=\small},
        legend style={font=\footnotesize, legend columns=4, legend cell align=center},
        label style={font=\small, text width=4.0cm, align=center},
        width=6cm,}
      \inputaxis{\dirdefault/gaussian_processoutput_index_1_rb_size_63_dual_rb_size_21_v1_tikz.tex}{%
        name=gpviz1, legend to name=leg:gpviz1big,%
        title={}}
      \inputaxis{res2bias/gaussian_processoutput_index_1_rb_size_63_dual_rb_size_21_v1_tikz.tex}{%
        at={($(gpviz1.east)+(3cm,0)$)}, anchor=west, legend to name=leg:gpviz2big,%
        title={}}
      \node at ($(gpviz1.south east)+(1.5cm, -1cm)$) [anchor=north] {\pgfplotslegendfromname{leg:gpviz1big}};
    }
    \end{tikzpicture}

    \caption{
    Visualization of \methodname{} surrogates ($\error = \enstateerror{}$ and
    $\outbias{}$, $\indicators = \log{\residualnorm}$, $\functionErrorNo =
    \log$), computed using $\romesTraining =100$ training points and the GP kernel method for
    a higher dimensional ROM with $p=62$.}\label{fig:gpvisualizationBig}
  \end{center}
\end{figure}

  \begin{figure}[htp]
    \begin{center}
      \tikzsetfigurename{effectivityAndImprovementBig}
      \begin{tikzpicture}
      {
        \pgfplotsset{
          tick label style={font=\small},
          legend style={font=\small, legend columns=-1, legend cell align=center},
          label style={font=\small, text width=4.5cm, align=center},
          width=6cm}
        \inputaxis{\dirdefault/effectivity_mode_50_bs_\rblarge_v\vardefault_tikz.tex}{%
          legend to name=leg:effectivityBig, name=effectivity1}
        \inputaxis{res2bias/improvement_mode_bs_\rblarge_v\vardefault_tikz.tex}{%
          at={($(effectivity1.east)+(3cm,0)$)}, anchor=west,%
          restrict y to domain=-2:20,%
          ymin=0, ymax=3,%
          legend to name=leg:improvementBig}
        \node at ($(effectivity1.south east)+(1.5cm, -1cm)$) [anchor=north] {\pgfplotslegendfromname{leg:effectivityBig}};
      }
      \end{tikzpicture}

    \end{center}
    \caption{Effectivity $\eta(0.5,\bmu)$ of \methodname{} surrogate 
    ($\error = \enstateerror{}$, $\indicators = \log{\residualnorm}$, $\functionErrorNo =
\log$) and expected improvement
		$I(\stochastic{\delta}, \bmu)$ of \methodname{} surrogate
    (GP
		kernel, compliant $\delta = \outbias$, $\indicators = \log{\residualnorm}$, $\functionErrorNo =
\log$)  for a higher dimensional ROM with $\nrb = 62$
    and a varying number of training points
		$\romesTraining $.}\label{fig:effectivityAndImprovementBig}
  \end{figure}
\Figure{\ref{fig:effectivityAndImprovementBig}} reports the performance of these
	surrogates.
  Comparing the leftmost plot of 
	\Figure{\ref{fig:effectivityAndImprovementBig}} with
	\Figure{\ref{fig:effectivityComparisonRVM}}(i) shows that the state-space error surrogate
	exhibits nearly identical convergence for the larger- and smaller-dimension reduced-order models.
	As in the
	experiments of \Section{\ref{sec:experiments:comparison:RBbounds}}, the 
	value of 
	$\mode{\eta(0.5,\bmu)}$
  is close to $1$ in the mean. Comparing the rightmost plot of
	\Figure{\ref{fig:effectivityAndImprovementBig}} with
	\Figure{\ref{fig:improvementComparisonEldred}}(i) shows that the
	expected improvement for the output-error
  correction with the surrogate $\stochastic{\outbias}$ is around $0.2$ in the
  mean for both the larger- and smaller-dimension reduced-order models. However,
	for the larger-dimension reduced-order model,
  more training points are required to reduce the occurrence of
  improvement factors larger than $1$. This is an artifact of the low errors
	already produced by the larger-dimension ROM itself (i.e., small denominator in Eq.~\eqref{eq:improvementValueMF}).

	We therefore conclude that the \methodname{} method is applicable to ROMs of
	different fidelity.

\newcommand{\resdir}{dual_multi_output_rvm_no_transformation}
\newcommand{\resdirLog}{dual_multi_output_gp_log_transformation}
\newcommand{\outputOne}{output_index_3_rb_size_11_dual_rb_size_20}
\newcommand{\outputOneSmallRB}{output_index_3_rb_size_11_dual_rb_size_10}
\newcommand{\outputOneMediumRB}{output_index_3_rb_size_11_dual_rb_size_15}
\newcommand{\outputTwoMediumRB}{output_index_6_rb_size_11_dual_rb_size_17}
\newcommand{\outputTwoSmallRB}{output_index_6_rb_size_11_dual_rb_size_11}
\newcommand{\outputTwo}{output_index_6_rb_size_11_dual_rb_size_23}

\subsection{Multiple and non-compliant outputs}\label{sec:experiments:multipleOutputs}

Finally, we assess the performance of \methodname{} on a model with multiple
and non-compliant output functionals as discussed in
\Section{\ref{sec:dualWeightedResiduals}}. For this
experiment, we set two outputs to be temperate measurements at
points $x_1$ and $x_2$:
\begin{gather}
  \label{eq:multipleOutputFunctional}
  \out_i\rb{\bmu} \defeq \outputfun_i\rb{\state\rb{\bmu}} \defeq
	\outputfunCont_i\left(\contstate\rb{\bmu}\right) = \int_\Omega\dirac(x -
	x_i)\contstate\rb{x_i; \bmu}\text{d}x = \contstate\rb{x_i; \bmu}, \quad i = 1,2.
\end{gather}
where $\dirac$ denotes the Dirac delta function.
In this case, we construct a separate \methodname{} surrogates for each
output error $\stochastic{{\outbiasone}}$ and
$\stochastic{\outbiastwo}$. As previously discussed,
we use dual-weighted residuals as
indicators $\indicators_i(\bmu) = \redadjointcoarsetofineII(\bmu)^t \modresidual\left(
\redstate;\bmu\right)$, $i=1,2$ and
no transformation $d\equiv \text{id}_{\RR}$.
This necessitates the computation of approximate dual solutions, for which
dual reduced-basis spaces must be generated in the offline stage.  The
corresponding finite element problem can be found in Eq.~\eqref{eq:weakDual},
where Eq.~\eqref{eq:multipleOutputFunctional} above provides the right-hand
sides. The algebraic problems can be inferred from
Eq.~\eqref{eq:discreteDual}, where the discrete right-hand sides
are canonical unit vectors because the points $x_1$ and
$x_2$ coincide with nodes of the finite-element mesh.

Like the primal reduced basis, the dual counterpart can be generated with a
greedy algorithm that minimizes the approximation error for the reduced dual
solutions.

To assess the ability for \emph{uncertainty control} with the dual-weighted-residual
indicators (see Remark \ref{rem:uncertaintyControlDual}) we generate three
dual reduced bases of increasing fidelity: 1) error tolerance of $1$
(basis sizes $\nrbdual$ of 10 and
11), 2) error tolerance of $0.5$ (basis sizes $\nrbdual$ of 15 and
17), 3) error tolerance of $0.1$ (basis sizes $\nrbdual$ of 20 and
23).

To train the surrogates, we compute 
$\outbiasone{(\bmu)}$, $\outbiastwo{(\bmu)}$,
$\indicators_1(\bmu)$ (of varying fidelity), $\indicators_2(\bmu)$ (of varying fidelity),   for  $\bmu\in\Pararomes\subset\Para$ with
  $\card{\Pararomes} = 500$. The first $T=100$
  points define the training set $\Paratrainromes\subset\Pararomes$  and the
	following $400$ points constitute the validation
  set $\Paravalidate\subset\Pararomes$.

Figure
\ref{fig:correlationDual} depicts the observed relationship between indicators $\indicators_1(\bmu)$ (of
different fidelity) and the 
error in the first output $\outbiasone(\bmu)$. Note that as the dual-basis
size $\nrbdual$ increases, the output error exhibits a nearly exact linear dependence on
the dual-weighted residuals.  This is expected, as the residual operator is linear in the
state.
Therefore, the RVM with a linear polynomial basis
produces the best (i.e., minimum variance) results for the \methodname{}
surrogates in this case.
\begin{figure}[t]
  \begin{center}
    \tikzsetfigurename{correlation_dual}
    \begin{tikzpicture}
      {
      \pgfplotsset{
      tick label style={font=\small},
      label style={font=\small, text width=4.0cm, align=center},
      title style={font=\small},
      legend style={font=\small},
      every x tick scale label/.style={
          at={(1,0)},xshift=1pt,anchor=south west,inner sep=0pt
      },
      width=4.5cm,
      }

      \inputaxis{\resdir/correlation_1d\outputOneSmallRB_tikz.tex}{%
      name=dualcorrelation_small,%
      title={(i) dual RB size $p_y=10$}}

      \inputaxis{\resdir/correlation_1d\outputOneMediumRB_tikz.tex}{%
      at={($(dualcorrelation_small.east)+(2.2cm,0)$)}, anchor=west,%
      name=dualcorrelation_middle,%
      title={(ii) dual RB size $p_y=15$}}

      \inputaxis{\resdir/correlation_1d\outputOne_tikz.tex}{%
      at={($(dualcorrelation_middle.east)+(2.2cm,0)$)}, anchor=west,%
      name=dualcorrelation,%
      title={(iii) dual RB size=$p_y=20$}}

      }
    \end{tikzpicture}
  \end{center}
  \caption{Relationship between between dual-weighted-residual indicators
	$\indicators_1 = \redadjointcoarsetofinearg{1}(\bmu)^t \modresidual\left(
\redstate;\bmu\right)$
	and errors in the (non-compliant) first output $\outbiasone$.
  \label{fig:correlationDual}}
\end{figure}

\Figure{\ref{fig:dualImprovement}} reflects the necessity of employing a
large enough dual reduced basis to compute the dual-weighted-residual error
indicators.  For a small dual reduced basis, there is almost no improvement in
the mean, and only a slight improvement in the median; in some cases, the
`corrected' outputs are actually less accurate. However, the most accurate
dual solutions yield a mean and median error improvement of two orders of
magnitude.
	This illustrates the ability
and utility of
\emph{uncertainty control} when dual-weighted residuals are used as error
indicators.

\begin{figure}[thp]
  \begin{center}
    \tikzsetfigurename{improvementDual}
    \begin{tikzpicture}
      {
      \pgfplotsset{
        tick label style={font=\small},
        legend style={font=\small, legend columns=-1, legend cell align=center},
        title style={font=\small},
        label style={font=\small, text width=4.5cm, align=center},
        xmax=80,
        width=4.7cm}

      \inputaxis{\resdir/improvement_OutputOneSmall_tikz.tex}{%
      name=dualImprovementOneSmall,%
      legend to name=leg:dualImprovementOneSmall,%
      ymode={log},
      ylabel={improvement (smaller=better)},%
      ymin=1e-3, ymax=1e1,
      title={(i) dual RB size $p_y=10$}}

      \inputaxis{\resdir/improvement_OutputOneMedium_tikz.tex}{%
      at={($(dualImprovementOneSmall.east)+(1.7cm,0)$)}, anchor=west,%
      name=dualImprovementOneMedium,%
      legend to name=leg:dualImprovementOneMedium,%
      ymode={log},
      ylabel={},
      ymin=1e-3, ymax=1e1,
      title={(ii) dual RB size $p_y=15$}}

      \inputaxis{\resdir/improvement_OutputOneLarge_tikz.tex}{%
      at={($(dualImprovementOneMedium.east)+(1.7cm,0)$)}, anchor=west,%
      name=dualImprovementOne,%
      legend to name=leg:dualImprovementOne,%
      ymode={log},
      ylabel={},
      ymin=1e-3, ymax=1e1,
      title={(iii) dual RB size $p_y=20$}}

      \inputaxis{\resdir/improvement_OutputTwoSmall_tikz.tex}{%
      at={($(dualImprovementOneSmall.south)+(0,-2cm)$)}, anchor=north,%
      name=dualImprovementTwoSmall,%
      legend to name=leg:dualImprovementTwoSmall,%
      ymode={log},
      ylabel={improvement (smaller=better)},%
      ymin=1e-3, ymax=1e1,
      title={(iv) dual RB size $p_y=11$}}

      \inputaxis{\resdir/improvement_OutputTwoMedium_tikz.tex}{%
      at={($(dualImprovementTwoSmall.east)+(1.7cm,0)$)}, anchor=west,%
      name=dualImprovementTwoMedium,%
      legend to name=leg:dualImprovementTwoMedium,%
      ymode={log},
      ylabel={},
      ymin=1e-3, ymax=1e1,
      title={(v) dual RB size $p_y=17$}}

      \inputaxis{\resdir/improvement_OutputTwoLarge_tikz.tex}{%
      at={($(dualImprovementTwoMedium.east)+(1.7cm,0)$)}, anchor=west,%
      name=dualImprovementTwo,%
      legend to name=leg:dualImprovementTwo,%
      ymode={log},
      ylabel={},
      ymin=1e-3, ymax=1e1,
      title={(vi) dual RB size $p_y=23$}}
    \node at ($(dualImprovementTwoMedium.south)+(0, -1.5cm)$) {%
    \pgfplotslegendfromname{leg:dualImprovementTwo}};

      }
    \end{tikzpicture}
  \end{center}
  \caption{
		Expected improvement $I(\stochastic{\delta}, \bmu)$ for \methodname{}
		surrogate
		 (RVM, $\delta = \outbias$, $\indicators_i = \redadjointcoarsetofineII(\bmu)^t \modresidual\left(
\redstate;\bmu\right)$, $i=1,2$,
		$\functionErrorNo = \text{id}_{\RR}$)
		for a varying
		number of training points $T$
		and different dual reduced-basis-space dimensions. Compare with
		\Figure{\ref{fig:improvementComparisonEldred}} (1: no improvement, $>1$: error
  worsened, $<1$: error improved).
  \label{fig:dualImprovement}}
\end{figure}

Table
\ref{tab:dualValidation} reports validation results for the inferred
confidence intervals. While the validation results are quite good (and
appear to be converging to the correct values), they are not as accurate as those obtained for the compliant
output.  

\begin{figure}[thp]
  \begin{center}
    \tikzsetfigurename{ConfidenceValidationDual}

    \begin{tikzpicture}
      {
      \node (c11) {
      \small
       \input{\imagefile{\resdir/confidence_interval_table_tikz.tex}}
       };
       }
    \end{tikzpicture}
  \end{center}
  \caption{Gaussian-process validation for the \methodname{} surrogates (RVM, $\delta = \outbias$, $\indicators_i = \redadjointcoarsetofineII(\bmu)^t \modresidual\left(
\redstate;\bmu\right)$, $i=1,2$,
		$\functionErrorNo = \text{id}_{\RR}$).
    The table reports how often the actual error lies in the inferred confidence
    intervals. 
	The largest dual reduced-basis space
  dimensions ($\nrbdual = 20$ and $\nrbdual = 23$) are used to compute the
	error indicators. \label{tab:dualValidation}. }
\end{figure}

\section{Conclusions and outlook}

	This work presented the \methodname{} method for statistically modeling
	reduced-order-model errors. In contrast to rigorous error bounds, such
	statistical models are useful for tasks in uncertainty quantification.
	The method employs supervised machine learning methods to
	construct a mapping from existing, cheaply computable ROM error indicators to a
	\emph{distribution} over the true error. This distribution reflects the
	epistemic uncertainty introduced by the ROM.  We proposed \methodname{}
	ingredients (supervised-learning method, error indicators, and transformation
	function) that yield low-variance, numerically validated models for
	different types of ROM errors.
	
	For normed outputs, the
	\methodname{} surrogates led to effectivities with low variance and means close
	to the optimal value of one, as well as a notion of probabilistic rigor.
	This is in contrast to existing ROM error bounds, which exhibited mean
	effectivities close to ten; this improvement will likely be more pronounced
	for more complex (e.g., nonlinear, time dependent) problems. Further, the
		\methodname{} surrogates were computationally less expensive than the
		error bounds, as the
		coercivity-constant lower bound was not required.

For general outputs, the \methodname{} surrogate allowed the ROM output to be
corrected, which yielded a near 10x accuracy improvement. Further, the
uncertainty in this error could be \emph{controlled} by modifying the
dimension of the dual reduced basis. On the other hand, existing approaches
(i.e., multifidelity correction) that employ system inputs (not error
indicators) as inputs to the error model did not lead to improved output
predictions. This demonstrated the ability of \methodname{} to mitigate the
curse of dimensionality: although the problem was characterized by nine system
inputs, only one error indicator was necessary to construct a low-variance,
validated \methodname{}
surrogate.
	We foresee the combination of ROMs with \methodname{} error
	surrogates to be powerful in UQ applications, especially
	when the number of system inputs is large. 
  Future work entails integrating and analyzing the 
  ROM/\methodname{} combination for specific UQ problems, e.g., Bayesian
	inference, as well as automating the procedure for 
selecting \methodname{} ingredients for different problems.
	Future work will also involve integrating the \methodname{} surrogates into
	the greedy algorithm for the reduced-basis-space selection; this has the
	potential to improve ROM quality due to the near-optimal expected
	effectivities of the error surrogates.

\section*{Acknowledgments}
We thank Khachik Sargsyan for his support in the understanding, selection, and
development of the machine-learning algorithms.  
This research was supported in part by an appointment to the Sandia National
Laboratories Truman Fellowship in National Security Science and Engineering,
sponsored by Sandia Corporation (a wholly owned subsidiary of Lockheed Martin
Corporation) as Operator of Sandia National Laboratories under its U.S.
Department of Energy Contract No. DE-AC04-94AL85000.  We also
acknowledge support by the Department of Energy Office of Advanced
Scientific Computing Research under contract 10-014804.

\ifthenelse{\equal{\atype}{preprint}}
{\appendix

\section{Reduced-basis method for parametric elliptic PDEs with affine
parameter dependence}
\label{sec:rb}
This section summarizes the reduced-basis method for generating a
reduced-order model for parametric elliptic PDEs with affine parameter
dependence.
See Ref.~\cite{PR07} for additional details.

\newcommand{\InnerProductM}{\mathbf{K}}
\newcommand{\Riesz}{\InnerProductM^{-1}}
\newcommand{\contrbstate}{u_{\redI}}

\subsection{Parametric elliptic problem}
\label{sec:rb:ellipticProblem}
We start with the discretized version of the elliptic problem, where
the state resides in a function space $X_{\detI}\subset X$ of finite
dimension $\ndof\defeq \dim(X_{\detI})$, e.g., a finite-element space. For
details on the derivation of such a finite element discretization from the
analytical formulation of an elliptic PDE and its convergence properties, we
refer to Ref.~\cite{Braess2007}.

The problem reads as follows: For any point in the input space $\bmu
\in \Para$, compute a quantity of interest
\begin{align}
  \outputs(\bmu)&\defeq\outputfunCont(\contfestate(\bmu))
  \label{eq:ellipticFOMOutput} \\
  \intertext{with $\contfestate\in X_{\detI}$ fulfilling}
  a(\contfestate(\bmu), \contfetest;\bmu) &= f(\contfetest;\bmu) \ \text{for
	all }\contfetest\in X_{\detI}.
  \label{eq:ellipticFOMEquation}
\end{align}
Here, for all inputs $\bmu \in \Para$, $a(\cdot, \cdot;
\bmu):X_{\detI}\times X_{\detI}\to X_{\detI}$ is a bilinear form on a Hilbert
space $X_{\detI}$ and is
symmetric, continuous, and coercive, fulfilling
\begin{equation}
  \sup_{u\in X} \sup_{v \in X} \frac{a(u,v;\bmu)}{\norm{u}_{X}\norm{v}_{X}} < \gamma(\bmu) < \infty
  \quad
  \text{ and }
  \inf_{u\in X, \norm{u}_X=1} a(u, u;\bmu) > \alpha(\bmu)
\end{equation}
for constants $\gamma(\bmu)$
and $\alpha(\bmu) > 0$. In addition, we have $f(\bmu) \in X^*_{\detI}$.

Given a basis $\cb{\febasisvec_i}_{i=1}^\ndof$ spanning $X_{\detI}$, the
solution can be expressed as a linear combination $\contfestate(\bmu) = \sum_{i=1}^\ndof
{u}_i(\bmu) \febasisvec_i\in X_{\detI}$ with the state
vector $\state\defeq \left[u_1 \ \cdots\ u_\ndof\right]^t :\Para\rightarrow \RR^\ndof$
containing the solution degrees of freedom. Then, problem \eqref{eq:ellipticFOMEquation} can be solved
algebraically by solving the system of linear equations 
\begin{equation}
  \residual(\state;\bmu)\defeq \mathbf{A}(\bmu) \state(\bmu) - \mathbf{f} (\bmu) = 0,
  \label{eq:ellipticFOMEquationLES}
\end{equation}
with $\mathbf{A}:\Para\rightarrow \RR^{\ndof\times \ndof}$ given by $\mathbf{A}_{ij}(\bmu) =
a\rb{\febasisvec_i, \febasisvec_j;\bmu}$ and $\mathbf{f}:\Para\rightarrow \RR^\ndof$ given by
$\mathbf{f}_i(\bmu) = f(\febasisvec_i;\bmu)$. Similarly, the output can be expressed as a
function of the state degrees of freedom:
\begin{equation} 
\outputfun(\state(\bmu)) \defeq\outputfunCont\left(\sum_{i=1}^\ndof
{u}_i(\bmu) \febasisvec_i\right).
\end{equation}

The key concept underlying projection-based model-reduction techniques is to
find a lower dimensional subspace $X_{\redI} \subset X_{\detI}$ of dimension
$\nrb: = \dim(X_{\redI}) \ll \ndof$, such that \eqref{eq:ellipticFOMEquationLES} can
be solved more efficiently.
For the reduced-basis method, this problem-dependent {\em reduced-basis
space} is obtained through identification of a few training points
$\Paratrain\defeq\cb{\bmutrain_i}_{i=1}^p\subset\Para$ for which the
full-order equations \eqref{eq:ellipticFOMEquationLES} are solved.
These functions then span the reduced basis space
$X_{\redI} = \myspan\cb{u(\bmutrain_i)}_{i=1}^p$
from which we generate an orthonormal basis
$\cb{\rbphi_i}_{i=1}^{\nrb}.$ 
We refer to subsection \ref{sec:rb:greedy} for more details on the selection
process of these parameters. Then, these basis vectors can be expressed in
terms of the original finite-element degrees of freedom, i.e.,
\begin{equation} 
\rbphi_{j} = \sum_{i=1}^\ndof\febasisvec_i\rightbasis_{ij}, \quad
j=1,\ldots,\nrb.
\end{equation} 
This transformation provides the entries of the trial basis $\rightbasis$. As
we consider only Galerkin projection (due to the symmetric and coercive nature of the
PDE), we set $\leftbasis=\rightbasis$, which leads to a reduced-order problem
of the form
\eqref{eq:ROMdyn}:
\begin{equation}
\rightbasis^t\residual(\rightbasis\redstatecoords;\bmu) =
\rightbasis^t\mathbf{A}(\bmu)\rightbasis\redstatecoords(\bmu) -
\rightbasis^t\mathbf{f}(\bmu)\rightbasis = 0.
  \label{eq:ellipticROMEquationLES}
\end{equation}
The reduced quantity of interest is then given by $\redout{} = \outputfun(\rightbasis \redstatecoords).$
\subsection{Offline--online decomposition}\label{sec:offlineOnline}

The low-dimensional operators in Eq.~\eqref{eq:ellipticROMEquationLES} can be
computed efficiently assuming that the bilinear forms $a$ and the functional
$f$ are {\em affinely parameter dependent}, i.e., they can be written as
\begin{equation}
  \label{eq:affinelyParameterDependentCondition}
  a(\cdot, \cdot; \bmu) = \sum_{q=1}^{Q_a} \sigma_a^q(\bmu) a^q(\cdot, \cdot) \quad\text{and}\quad
  f(\cdot; \bmu) = \sum_{q=1}^{Q_f} \sigma_f^q(\bmu) f^q(\cdot).
\end{equation}
with parameter-independent bilinear forms $a^q:X_h \times X_h \to X_h,
q=1,\ldots,Q_a$ and functionals $f^q:X_h\to X_h, q=1,\ldots, Q_f$, and
coefficient functions $\sigma_a^q:\paramDomain \to \RR, q=1,\ldots,Q_a$ and
$\sigma_f^q:\paramDomain \to \RR, q=1,\ldots, Q_f$.

In this affinely parameter dependent case, the parameter-dependent reduced
quantities
$\rightbasis^t\mathbf{A}(\bmu)\rightbasis$ and
$\rightbasis^t\mathbf{f}(\bmu)\rightbasis$ can be quickly assembled via a linear
combination of their parameter-independent components:
\begin{gather}
    \rightbasis^t\mathbf{A}(\bmu)\rightbasis
    =
    \sum_{q=1}^{Q_a} \sigma^q_a(\bmu) \rightbasis^t\mathbf{A}^q\rightbasis\\
    \rightbasis^t\mathbf{f}\rightbasis
    =
    \sum_{q=1}^{Q_f} \sigma^q_f(\bmu) \rightbasis^t\mathbf{f}^q\rightbasis.
\end{gather}

Here the matrices $\mathbf{A}^q \in \RR^{\ndof\times\ndof}, q=1,\ldots,Q_a$ are
given by $\mathbf{A}^q_{ij} = a^q(\febasisvec_{i}, \febasisvec_{j})$ and the
vectors $\mathbf{f}^q \in \RR^{\ndof}, q=1,\ldots,Q_f$ are given by
$\mathbf{f}^q_{i} = f(\febasisvec_{i})$.

The parameter-independent components
\begin{equation}
  \rightbasis^t\mathbf{A}^q\rightbasis \in \RR^{p \times p},\ q=1,\ldots,Q_a
  \qquad \text{and} \qquad
  \rightbasis^t\mathbf{f}^q\rightbasis \in \RR^{p},\ q=1,\ldots,Q_f
\end{equation}
must be computed only once, during the so-called {\em
offline stage}. All parameter-dependent
computations are carried out in the {\em online stage} in an efficient manner,
i.e., with computational complexities that only depend on the reduced-basis
dimension $p$. We refer to Ref.~\cite{PR07} for more details on the
offline/online decomposition. 
Note that in the non-affine case, low-dimensional operators can be
approximated using the {\em empirical interpolation method} \cite{BMNP04, HO08c, CS09}.

\subsection{Error bounds}
\label{sec:ellipticProblem:errorBounds}

The state error can be measured either in the problem--dependent
energy norm or in the norm given by the underlying Hilbert space $X_{\detI}$.
The energy norm for the above problem is defined as
\begin{equation}
  \enorm{u} \defeq a(u,u;\params)
\end{equation}
and is equivalent to the $X_h$-norm 
\begin{equation}
  \alpha \xnorm{u}^2 \leq \enorm{u}^2 \leq \gamma \xnorm{u}^2.
\end{equation}
The bounds for the state errors $\statebound^{\bmu}$ and $\statebound$ for the errors
$\enstateerror{}$ and $\xstateerror{}$ are defined as follows:
\begin{alignat}{2}
  \statebound^{\bmu}(\bmu) &\defeq \frac{r(\bmu)}{\sqrt{\alpha_{\text{LB}}(\bmu)}} & &\geq \enstateerror{(\bmu)} \text{ and}
  \label{eq:stateboundDefinitionEnergy}
  \\
  \statebound{(\bmu)} &\defeq \frac{r(\bmu)}{\alpha_{\text{LB}}(\bmu)}               & &\geq \xstateerror{(\bmu)},
  \label{eq:stateboundDefinition}
\end{alignat}
where $\alpha_{\text{LB}}(\bmu) \leq \alpha(\bmu)$ is a lower bound for the
coercivity constant. This lower bound can easily computed if the bilinear form
$a$ is symmetric and {\em parametrically coercive} as defined in
\cite[Ch.4.2]{PR07}, i.e., there is a point $\bar{\bmu} \in \paramDomain$,
such that
\begin{gather}
  (u, v)_{X_h} = a(u, v; \bar{\bmu}), \\
  \intertext{the parameter independent bilinear forms are coercive}
  a^q(w, w) \geq 0 \qquad \forall \bmu \in \paramDomain, q=1, \ldots, Q_a \\
  \intertext{and their coefficients strictly positive}
  \sigma_a^q(\bmu) > 0 \qquad \forall \bmu\in \paramDomain, q=1,\ldots, Q_a.
\end{gather}

In this special case, the lower coercivity constant can be chosen as
\begin{equation}
  \alpha_{\text{LB}}
  = \min_{q=1,\ldots,Q_a} \frac{\sigma^q_a(\bmu)}{\sigma^q_a(\bar{\bmu})}.
  \label{eq:mincoercivityAlphaLB}
\end{equation}

A similar bound can be generated for errors in \textit{compliant} output
functionals $\outputfunCont=f$ (and $a$ symmetric):
\begin{equation}
  \outbound{(\bmu)} \defeq \frac{r(\bmu)^2}{\alpha_{\text{LB}}(\bmu)} \geq
  \outbias{(\bmu)} \defeq \out(\bmu) - \redout{(\bmu)} > 0.
  \label{eq:outboundDefinition}
\end{equation}

Analogously to the upper bounds, lower bounds for the errors can be computed as

\begin{alignat}{2}
  \statebound^{LB, \bmu}(\bmu) &\defeq \frac{r(\bmu)}{\sqrt{\gamma_{\text{UB}}(\bmu)}} & &\leq \enstateerror{(\bmu)},
  \label{eq:statelowerboundDefinitionEnergy}
  \\
  \statelowerbound{(\bmu)} &\defeq \frac{r(\bmu)}{\gamma_{\text{UB}}(\bmu)}               & &\leq \xstateerror{(\bmu)} \text{ and}
  \label{eq:statelowerboundDefinition}
  \\
  \outlowerbound{(\bmu)} &\defeq \frac{r(\bmu)^2}{\gamma_{\text{UB}}(\bmu)} &&\leq
  \outbias{(\bmu)}.
\end{alignat}

where $\gamma_{\text{UB}} \geq \gamma$ is an upper bound for the continuity
constant $\gamma$.

\subsubsection{Offline--online decomposition of the residual norm}\label{sec:offlineOnlineResidual}

In order to compute the error bounds efficiently, the residual norm
\begin{equation}
  r(\bmu) \defeq \norm{f(\cdot; \bmu) - a(\contrbstate(\bmu), \cdot; \bmu)}_{X_h^*}
\end{equation}
must be computed efficiently through an offline/online decomposition. Here, the
reduced solution $\contrbstate(\bmu)$ is given by $\sum_{i=0}^{\nrb} \hat{u}_i
\rbphi_i \in X_{\redI} \subset X_{\detI}$ with the reduced state vector
$\redstatecoords = \left[\hat u_1 \ \cdots\ u_\nrb\right]^t :
\Para\rightarrow \RR^\nrb$ containing the reduced solution degrees of
freedom. Algebraically the residual norm can be computed as
\begin{equation}
  \begin{gathered}
    \rb{r(\bmu)}^2
    =
    (\Riesz \residual{(\rightbasis\redstatecoords(\params);\params}))^t
     \residual{(\rightbasis\redstatecoords(\params);\params})
     \\
     = \sum_{1\leq q,q' \leq Q_a} \sigma^q_a(\bmu)\sigma^{q'}_a(\bmu)\;\;
     \redstatecoords(\params)
     \rb{\rb{\Riesz \mathbf{A}^q \rightbasis}^t \mathbf{A}^{q'} \rightbasis}
     \redstatecoords(\params)
     \\
     - 2 \sum_{q=1}^{Q_f} \sum_{q'=1}^{Q_a} \sigma^q_f(\bmu)\sigma^{q'}_a(\bmu)\;\;
     \rb{\rb{\Riesz \mathbf{f}^q \rightbasis}^t \mathbf{A}^{q'} \rightbasis}
     \redstatecoords(\params)
     \\
     + \sum_{1\leq q,q' \leq Q_f} \sigma^q_f(\bmu)\sigma^{q'}_f(\bmu)\;\;
     + \rb{\rb{\Riesz \mathbf{f}^q \rightbasis}^t \mathbf{f}^{q'} \rightbasis}
\end{gathered}
\end{equation}
where the inner product matrix $\InnerProductM \in \RR^{\ndof \times \ndof}$ is
given by $\InnerProductM_{ij} = (\febasisvec_i, \febasisvec_j)_{X_h}$. It
is used to identify the Riesz representation of the residual $f(\cdot;
\bmu) - a(\contrbstate(\bmu), \cdot; \bmu) \in X_\detI^*$.

All low-dimensional and parameter-independent matrices
\begin{gather}
  \rb{\rb{\Riesz \mathbf{A}^q \rightbasis}^t \mathbf{A}^{q'} \rightbasis}
  \in \RR^{\nrb\times\nrb}, \qquad 1\leq q,q' \leq Q_a\\
  \rb{\rb{\Riesz \mathbf{f}^q \rightbasis}^t \mathbf{A}^{q'} \rightbasis}
  \in \RR^{\nrb}, \qquad 1\leq q,q' \leq Q_f,Q_a\\
  \rb{\rb{\Riesz \mathbf{f}^q \rightbasis}^t \mathbf{f}^{q'} \rightbasis}
  \in \RR, \qquad 1\leq q, q' \leq Q_f
\end{gather}
can be pre-computed during the {\em\/ offline stage}.

\subsubsection{Dual weighted error estimates}\label{sec:supp:dual}
As discussed in \Section{\ref{sec:dualWeightedResiduals}}, the \methodname{}
surrogate for
modeling errors in general (non-compliant) outputs $\outputfunCont$ requires a dual
solution $\dualstate(\params) \in \RR^{\ndof}$ . In the present context, the
associated problem is: Find $\contfedual(\bmu)\in X_h$ fulfilling
\begin{equation}\label{eq:weakDual}
  a(\contfetest, \contfedual(\bmu); \bmu) = -\outputfunCont(\contfetest) \qquad
	\text{for all }\contfetest \in X_h.
\end{equation}
Analogously to the discussion for the primal problem, we require a
reduced-basis space $X_{\redI}^{\outputfunCont} \subset X_h$ for this dual
problem. Algebraically, this leads to a reduced-basis matrix
$\rightbasisdualcoarse \in \RR^{\ndof\times\nrbdual}$ associated for a
reduced dual solution $\redadjointcoarse \in \RR^{\nrbdual}$ such that
\begin{equation}\label{eq:discreteDual}
  \rightbasisdualcoarseT\residual_{\outputfun}(\rightbasisdualcoarse\redadjointcoarse(\params); \params) = 0, \qquad \text{with}\qquad \residual_{\outputfun}(\adjoint;\params) \defeq \mathbf{A}(\params)\adjoint + \mathbf{\outputfun}.
\end{equation}
Here, $\mathbf{\outputfun} \in \RR^{\ndof}$ is given by
$(\mathbf{\outputfun})_i = \outputfunCont(\febasisvec_i)$.
Then, this error can be bounded as
\begin{equation}
  \outbound{(\bmu)} \defeq \frac{\residualnorm(\params)\residualnorm_{\outputfun}(\params)}{\alpha_{\text{LB}}(\params)} \geq
  |\outbias{(\bmu)}|,
  \label{eq:outboundDefinitionDual}
\end{equation}
with the dual residual norm $\residualnorm_{\outputfun}(\params) \defeq
\norm{\residual_{\outputfun}(\rightbasisdualcoarse
\redadjointcoarse(\params);\params)}_2$.

\subsection{Basis generation with greedy algorithm}
\label{sec:rb:greedy}

We have so far withheld the mechanism by which the reduced basis spaces $X_{\redI}, X_{\redI}^g
\subset X_{\detI}$ are constructed. This section provides a brief
overview on this topic, and to ease exposition, we will focus on the
generation of the primal reduced-basis space only. As previously mentioned,
the reduced-basis space is constructed from $\nrb$ chosen input-space points
$\Paratrain\subset\Para$ with $\card{\Paratrain} = \nrb$. As the reduced-basis
solutions should provide accurate
approximations for all points $\bmu \in \Para$, the ultimate
goal is to find a reduced-basis space of low dimension $\nrb\ll\ndof$ that minimizes some
distance measure between itself and the manifold $\cM\defeq\cb{u(\bmu)\ |\
\bmu \in \Para}$. Two possible candidates for such a distance measure are the
maximum projection error
\begin{equation}
  \label{eq:manifoldDistanceProjection}
  \text{dist}_1(X_{\redI}, \cM) \defeq \max_{u \in \cM} \inf_{v \in X_{\redI}} \norm{u-v}
\end{equation}
and the maximum reduced state error
\begin{equation}
  \label{eq:manifoldDistanceRB}
  \text{dist}_2(X_{\redI}, \cM) \defeq \max_{\contstate(\bmu) \in \cM} \norm{\contstate(\bmu)-\contredstate(\bmu)}.
\end{equation}
The optimally achievable distance of a reduced-basis space from the manifold
$\cM$ is known as the Kolmogorov $N$-width and is given by
\begin{equation}
  d_N(\cM)\defeq\inf_{\tilde{X} \subset X, \dim(\tilde{X})=N} \text{dist}_1(\tilde{X}, \cM).
\end{equation}
While the Kolmogorov $N$-width is a theoretical measure and is only known for a
few simple manifolds, it is possible to construct a reduced-basis space
with a myopic {\em greedy algorithm} that often comes close to this optimal value.
The
greedy algorithm works as follows: First, the manifold of `interesting
solutions' is discretized by defining a finite subset
$\Mgreedy\defeq\cb{\contstate(\params)\,|\,\params \in \Paragreedy} \subset
\cM$,  with $\card{\Paragreedy} \geq \nrb$ for
which the distance measures \eqref{eq:manifoldDistanceProjection} and
\eqref{eq:manifoldDistanceRB} can
be computed in theory. However, as computing these distance measures exactly
requires knowledge of the full-order solution $\contstate$, the method
substitutes these
distance measures by error bounds $\max_{\bmu \in \Paragreedy}
\statebound{(\bmu)}$. This allows for the construction of a sequence of reduced-basis
spaces $X_1 \subset X_2 \subset \cdots \subset X_\nrb =: X_{\redI}$ by choosing
the first subspace arbitrarily, as
\begin{equation}
  X_1 =\myspan\cb{u(\bmu_1)}
\end{equation}
with a randomly chosen parameter $\bmu_1 \in \Para$.
Then, the following spaces can be chosen iteratively by computing the parameter that maximizes
the error bound
\begin{equation}
  \bmu^k_{\max} \defeq \arg \max_{\bmu \in \Paragreedy} \statebound^k{(\bmu)}
\end{equation}
and adding the corresponding full-order solution to the reduced basis spaces
\begin{equation}
  X_{k+1} = X_k \cup \myspan\cb{u(\bmu^k_{\max})},
\end{equation}
where $\statebound^k{(\bmu)}$ denotes the bound for model-reduction errors
$\norm{\contstate(\bmu) - \contredstate(\bmu)}$ computed using
reduced-basis space $X_k$.
This algorithm works well in practice and has recently been verified by
theoretical convergence results: Refs.~\cite{BC11, Ha11} proved
that the distances of these heuristically constructed reduced
spaces converge at a rate similar to that of the theoretical
Kolmogorov $N$-width if this $N$-width converges algebraically or exponentially
fast.

In practice, the performance of this greedy algorithm depends
strongly on both the effectivity and computational cost of the error
bounds, as less expensive error bounds allow for a larger number of
points in $\Mgreedy$. Therefore, we expect that employing the \methodname{}
surrogates in lieu of the rigorous bounds may lead to performance gains for the
greedy algorithm; this constitutes an interesting area of future investigation.

}{}

\bibliography{reduced_basis,uq}

\end{document}